\documentclass[a4paper, reqno]{amsart}

\RequirePackage[T1]{fontenc}
\usepackage[utf8]{inputenc}
\usepackage{lmodern}
\usepackage[english]{babel}
\usepackage[linkcolor=blue, citecolor=green, colorlinks=true]{hyperref}
\usepackage{geometry}
\usepackage{fullpage}
\usepackage{graphicx}
\usepackage{xcolor}
\usepackage{amsmath, amssymb}
\usepackage{titlesec}
\usepackage{mathptmx}
\usepackage{fancyhdr}
\usepackage{authoraftertitle}
\usepackage{titling}

\usepackage{chngcntr}
\numberwithin{equation}{section}

\usepackage{etoolbox}

\newcommand\datefootnote[1]{%
  \begingroup
  \renewcommand\thefootnote{}\footnote{#1}%
  \addtocounter{footnote}{-1}%
  \endgroup
}

\let\MakeTitle\maketitle
\renewcommand\maketitle{%
  \MakeTitle
  \thispagestyle{empty}%
  \fontsize{14pt}{16pt}\selectfont
}
\makeatletter
\def\@maketitle{
\clearpage
  \null
  \vskip 2em%
  \begin{center}%
    {\Large\sffamily\bfseries \MakeUppercase\@title  \par}%
    \vskip .5em%
    {\large\mdseries\hspace{0cm}
    \vskip 1em
      \begin{tabular}[t]{c}%
        \sffamily\hspace{0cm} \theauthor
      \end{tabular}\par \large}%
    \vskip 1em%
    \datefootnote{\today}
    \pagenumbering{arabic}
  \end{center}%
}
\makeatother

\titleformat{\section}
        [block]
        {\Large\bfseries\sffamily\filcenter}
        {\thesection. }
        {0em}
        {\MakeUppercase}
        [
        ]

 \titleformat{\subsection}
        [block]
        {\Large\scshape}
        {\thesubsection}
        {0.5em}
        {}
        [
        ]
\titleformat{\subsubsection}
    [block]
    {\large\scshape}
    {\thesubsubsection}
    {0.5em}
    {\textbullet{} }
    []
    
    \usepackage{titletoc}
        \titlecontents{section}[1.5em]
       {\addvspace{0em}}
       {\contentslabel{2em}}
       {\hspace*{-2.3em}}
       {\hfill\contentspage}
   
\renewenvironment{abstract}{\bigskip \begin{center}\begin{minipage}[c]{0.8\textwidth} \normalsize}{\end{minipage} \end{center} \bigskip}

\makeatletter
\let\@oddfoot\@empty
\let\@evenfoot\@empty
\makeatother

\newcommand{\R}{\mathbf{R}}
\newcommand{\C}{\mathbf{C}}
\newcommand{\A}{\mathbf{A}}
\newcommand{\Q}{\mathbf{Q}}
\newcommand{\Z}{\mathbf{Z}}

\renewcommand{\c}{C}
\newcommand{\p}{P}
\renewcommand{\d}{D}

\newcommand{\n}{N}

\newcommand{\GL}{\mathrm{GL}}
\renewcommand{\O}{\mathcal{O}}
\newcommand{\legendre}[2]{
\left( \frac{#1}{#2} \right)}

\newcounter{pnum}[section]

\newtheorem{prop}{Proposition}
\newtheorem{coro}{Corollary}
\newtheorem{lem}{Lemma}

\newtheorem{thm}{Theorem}
\newtheorem{hyp}{Hypothesis}

\title{Quadratic twists of central~values~for~GL(3)}
\author{Chan Ieong Kuan and Didier Lesesvre}
\date{\today}
\address{School of mathematics (Zhuhai) \newline
Zhuhai Campus, Sun Yat-Sen University \newline
Tangjiawan, Zhuhai, Guangdong, 519082, China (PRC)}
\email{kuanchi3@mail.sysu.edu.cn}
\email{didier@mail.sysu.edu.cn, lesesvre@math.cnrs.fr}

\begin{document}

\maketitle

\begin{abstract}
We prove that a cuspidal automorphic representation of $\GL(3)$ over any number field is determined by the quadratic twists of its central value. In the case $\pi$ is not a Gelbart-Jacquet lift, the result is conditional on the analytic behavior of a certain Euler product. We deduce the nonvanishing of infinitely many quadratic twists of central values. This generalizes a result of Chinta and Diaconu that was valid only over $\Q$ and explored only for Gelbart-Jacquet lifts. 
\end{abstract}

\setcounter{tocdepth}{1}
\normalsize
\tableofcontents

\setlength{\parskip}{0.5em}	

\normalsize

\section{Introduction}
\label{sec:introduction}

\subsection{Statement of the result}
\label{subsec:statement-result}

Given an automorphic object, its associated $L$-function is built in a way to be a generating function with good analytic properties \cite{farmer_analytic_2017}. Of considerable importance is its value at special points that appears to carry a lot of information \cite{iwaniec_perspectives_2000}. For instance, in the automorphic forms setting, the non-vanishing at the central point of the $L$-function attached to a Rankin-Selberg product $\pi_1 \times \pi_2$ is conjectured to be related to the non-vanishing of automorphic periods of arithmetic significance \cite{cogdell_nonvanishing_2005}. With a geometric flavor, given an elliptic curve, the Birch and Swinnerton-Dyer conjecture relates the order of vanishing of its associated $L$-function at the central point to the rank of the elliptic curve. A last example coming from algebraic number theory states that a special value of the $L$-function attached to some characters is related to important invariants of the underlying field extension, by means of the class number formula. 

Following this philosophy, we are naturally interested in the following question: does the central value of the $L$-function attached to an automorphic form determine it? Obviously the single central value is insufficient since two non-isomorphic automorphic forms can be rescaled to have matching central values. A finite family of twists of the central value is also an elusive bet, since there are infinitely many linearly independent automorphic forms. Therefore, the most natural question appears to be the existence of an infinite family of twists of the central value that determines an automorphic form up to isomorphism. 

The first result in this direction is due to \cite[Theorem A]{luo_determination_1997}  for holomorphic modular forms. Let $S_k(N)$ be the space of holomorphic cusp forms of level $N$, weight $k$ and trivial central character. They proved the following: 

\begin{thm}[Luo-Ramakrishnan]
\label{thm:luo-ramakrishnan}
Let $f$ and $f'$ be normalized newforms in $S_{k}(N)$ and $S_{k'}(N')$ respectively. Let $L(s, f)$ and $L(s, f')$ be their associated completed $L$-functions. If there is a nonzero constant $\kappa$ such that for every quadratic character $\chi$ of conductor prime to $NN'$,
\begin{equation}
L\left(\tfrac{1}{2}, f \otimes \chi\right) = \kappa \cdot L\left(\tfrac{1}{2}, f' \otimes \chi\right),
\end{equation}
then $k=k'$, $N=N'$ and $f=f'$.
\end{thm}

 Based on the machinery developed by Fisher and Friedberg  for function fields \cite{fisher_double_2004} and adapting the methods developed by \cite{chinta_determination_2005} in the automorphic representations setting, Li proved an analogous result \cite[Theorem 1.1]{li_determination_2007} for automorphic representations of $\GL(2)$ over general number fields. More precisely, let $X$ be a family of quadratic characters as defined in \cite[Section 1.2]{li_determination_2007}. Li proved that an automorphic representation of $\GL(2)$ is determined by the associated quadratic twists of its central value.

\begin{thm}[Li]
\label{thm:li}
Let $\pi$ and $\pi '$ be two self-contragredient cuspidal automorphic representations of $\GL(2)$ over a number field $F$ with trivial central character. Let $L(s, \pi)$ and $L(s, \pi')$ be their associated completed $L$-functions. Assume there is a character $\chi \in X$ such that the $\varepsilon$-factor $\varepsilon(\frac{1}{2}, \pi \otimes \chi)$ is nonzero. If there exists a nonzero constant $\kappa$ such that for every $\chi \in X$,
\begin{equation}
L\left(\tfrac{1}{2}, \pi \otimes \chi\right) = \kappa \cdot   L\left(\tfrac{1}{2}, \pi' \otimes \chi\right),
\end{equation}
\noindent then $\pi \simeq \pi'$.
\end{thm}

Let $Y$ be the set of quadratic Dirichlet characters. In a similar fashion, \cite[Theorem 1.1]{chinta_determination_2005} proved the following result for cuspidal automorphic representations of $\GL(3)$ over $\Q$ that are Gelbart-Jacquet lifts:


\begin{thm}[Chinta-Diaconu]
\label{thm:chinta-diaconu}
Let $\pi$ and $\pi '$ be two self-contragredient cuspidal automorphic representations of $\, \GL(3)$ over $\Q$ with trivial central character. Let $L(s, \pi)$ and $L(s, \pi')$ be their associated completed $L$-functions. Assume that $\pi$ is a Gelbart-Jacquet lift. If there exists a nonzero constant $\kappa$ such that for every $\chi \in Y$,
\begin{equation}
L\left(\tfrac{1}{2}, \pi \otimes \chi\right) = \kappa \cdot  L\left(\tfrac{1}{2}, \pi' \otimes \chi\right),
\end{equation}
\noindent then $\pi \simeq \pi'$. 
\end{thm}

\textit{Remarks.} We note the following regarding the last two theorems: 
\begin{itemize}
\item[(\textit{i})] These results still hold if we consider partial $L$-functions instead of completed automorphic $L$-functions. It is also possible to restrict the family of characters to those with conductors relatively prime to a fixed number. 
\item[(\textit{ii)}] Li proves a slightly stronger result in the $GL(2)$ setting, allowing general central characters. As we are mostly interested in generalizing Chinta and Diaconu's result to general number fields and to explore the case of non-Gelbart-Jacquet lifts,  we do not add these technical variations.
\end{itemize}

One of the aims of the paper of Chinta and Diaconu was to explain the critical difficulties they faced when using their methods to more general number fields. Their proof heavily relies on the machinery of multiple Dirichlet series as developed by Bump, Diaconu, Friedberg, Goldfeld and Hoffstein \cite{diaconu_multiple_2003, bump_sums_2004}. However, the lack of good bounds for automorphic coefficients as well as mean square of central values proved to be a major obstacle for going beyond the field $\Q$ of rational numbers. Building on a quadratic large sieve inequality on number fields due to Goldmakher and Louvel \cite{goldmakher_quadratic_2013}, generalizing a result of Heath-Brown \cite{heath-brown_mean_1995}, we are able to overcome these difficulties. This is also the opportunity to fill in some technical gaps in the literature. 

The other difficulty is to go beyond the Gelbart-Jacquet assumption and to understand precisely the needed properties required to make the strategy work. Being a Gelbart-Jacquet lift is in fact equivalent to being self-contragredient, as a consequence of \cite{ramakrishnan_exercise_2014}, so that both hypotheses in \cite{chinta_determination_2005} are in fact redundant. Recent bounds on coefficients and non-vanishing results for symmetric square L-functions hence allow us to free the argument of the Gelbart-Jacquet (and the self-contragredient) assumption except for the following technical point:
\begin{hyp}
\label{hyp}
Let $a_\pi(N)$ be the coefficients of the L-function $L(s, \pi)$. Let $\alpha_\p, \beta_\p, \gamma_\p$ be the local spectral parameters of $\pi$. The Euler product
\begin{equation}
\label{hypothesis}
\prod_\p \left( 1+a_\pi(\p) |\p|^{-2s} \right)\left( 1-\alpha_\p \beta_\p |\p|^{-2s} \right)\left( 1- \beta_\p \gamma_\p |\p|^{-2s} \right)\left( 1 - \gamma_\p \alpha_\p |\p|^{-2s} \right)
\end{equation}

has no pole at $s=\tfrac{1}{2}$. Here, the product is taken over prime ideals of the ring of integers of the base field.
\end{hyp}

\textit{Remark.} Hypothesis \ref{hyp} is satisfied for Gelbart-Jacquet lifts. To the extent of our knowledge, we cannot quantify how weaker this hypothesis is when compared to $
\pi$ being a Gelbart-Jacquet lift.

We now state the main result of this paper.


\begin{thm}
\label{thm:result}
Let $\pi$ and $\pi'$ be two cuspidal automorphic representations of $\GL(3)$ over a number field $F$ with trivial central character. Let $L(s, \pi)$ and $L(s, \pi')$ be their associated completed $L$-functions. Let $X$ be the set of quadratic Hecke characters as defined in Section \ref{subsec:quadratic-symbols}. If there is a nonzero constant $\kappa$ such that for every $\chi \in X$,
\begin{equation}
L\left(\tfrac{1}{2}, \pi \otimes \chi\right) = \kappa \cdot   L\left(\tfrac{1}{2}, \pi' \otimes \chi\right),
\end{equation}
\noindent then 
\begin{itemize}
\item if $\pi$ is a Gelbart-Jacquet lift, then $\pi \simeq \pi'$;
\item if $\pi$ is not a Gelbart-Jacquet lift, if Hypothesis \ref{hyp} holds and if there is an ideal ray class $E$ such that  $L(1/2, \pi \otimes \chi_E) \neq 0$, then $\pi \simeq \pi'$.
\end{itemize}
\end{thm}

We discuss some aspects of our result below:
\begin{itemize}
\item[(i)] An automorphic representation $\pi$ of $\GL(3)$ with trivial central character is a Gelbart-Jacquet lift of a cuspidal automorphic representation of $\GL(2)$ \cite{ASENS_1978_4_11_4_471_0} if and only if its symmetric square $L$-function $L(s, \pi, \mathrm{sym}^2)$ has a pole at $1$, as proven in \cite{grs}. The difference in analytic behavior sparks the case separation of lifts and non-lifts.
\item[(ii)] Munshi and Sengupta \cite{munshi_determination_2015} recently proved that a $\GL(3)$ Hecke-Maass form is entirely determined by the twists of its central values by Dirichlet characters.  However, their result is restricted to the base field $\Q$ and requires all primitive characters.
\item[(iii)] There is no hope that quadratic twists of the central value can determine an automorphic representation, even in the case of $\GL(2)$: there are cuspidal automorphic representations $\pi$ on $\GL(2)$ with $L(1/2, \pi \otimes \chi)=0$ for every quadratic character $\chi$. One such example is given in the work of Waldspurger \cite{waldspurger_sur_1985}: take $\pi$ to be the base change lift to an imaginary quadratic field of the representation associated to an holomorphic cusp form on $\mathrm{SL}(2, \Z)$. A non-vanishing assumption is hence necessary and enforced in the $\GL(2)$ case in \cite{li_determination_2007}. Indeed, he assumes that there is a quadratic character such that $\varepsilon(1/2, \pi \otimes \chi)$ is one, which implies that there are infinitely many nonzero twists $L(1/2, \pi \otimes \chi)$ of the central value by \cite{friedberg_nonvanishing_1995}.
\item[(iv)] One could wonder whether a similar result exists for cuspidal automorphic representations of higher ranks $\GL(n)$ for $n \geqslant 4$. Our proof is based on the meromorphic continuation of the double Dirichlet series attached to the quadratic twists of the central values, and is essentially reduced to a double computation of the residue at the point $(\frac{1}{2},1)$.  A general fact in the theory of multiple Dirichlet series is that the group of known functional equations for $\GL(n)$, with $n \geqslant 4$, is insufficient to get a meromorphic continuation around this point, which implies that methods in this direction cannot address the problem for higher ranks. See \cite{bump_sums_2004} for more details on this spectral barrier problem.
\end{itemize}

We are mostly concerned in this article with the adaptation and methodology required to conquer the difficulties encountered by Chinta and Diaconu. As such, we assume familiarity with their paper --- we, however, always provide precise references therein.

\subsection{Non-vanishing consequence}


Non-vanishing questions are often difficult. \cite{luo_nonvanishing_2005} proved that an $L$-function $L(s, \pi)$ attached to an automorphic $L$-function on $\GL(n)$ over $\Q$ has infinitely many nonzero twists of its central value $L(1/2, \pi \otimes \chi)$ where $\chi$ is a Dirichlet character. However, nothing is said about the order of the character and showing non-vanishing results for quadratic characters remained out of reach. Similar to the $\GL(2)$ case, it is illusory to expect a non-vanishing result for quadratic twists without further assumptions. We prove the following result in the case of non-Gelbart-Jacquet lift representations, see Section \ref{nonvanishing}.

\begin{thm}
\label{thm:result-nv}
Let $\pi$ be a cuspidal automorphic representation of $\GL(3)$ over $F$. Assume that $\pi$ is not a Gelbart-Jacquet lift and that Hypothesis~\ref{hyp} is satisfied. If there is an ideal ray class $E$ such that $L(1/2, \pi \otimes \chi_E) \neq 0$, then there are infinitely many squarefree integral ideals $D$ such that $L(1/2, \pi \otimes \chi_D) \neq 0$.
\end{thm}

\subsection{Outlook of the proof}
\label{subsec:outlook}

The necessary knowledge on the quadratic characters $\chi_D$ for general number fields, parametrized by integral ideals $D$, and the background on the relevant automorphic $L$-functions are summarized in Section \ref{sec:background}. The structure of the proof of Theorem \ref{thm:result} relies on Chinta and Diaconu's original paper \cite{chinta_determination_2005} that was in the setting of rational numbers, and on the theory of multiple Dirichlet series as developed in \cite{bump_sums_2004}. The double Dirichlet series associated to the $L$-functions $L(s, \pi \otimes \chi_D)$, for quadratic characters $\chi_D$, is defined by
\begin{equation}
\widetilde{Z}(s, w; \pi) := \sum_{\substack{D}} \frac{L(s, \pi \otimes \chi_D)}{|D|^w} = \sum_{D, N} \frac{a_\pi(N)\chi_D(N)}{|N|^s|D|^w},
\end{equation}

\noindent where the sums run over nonzero integral ideals of $F$, the norm for ideals of $F$ is denoted $|\cdot |$, and the $a_\pi(N)$ are the Fourier coefficients of $\pi$. These double Dirichlet series are introduced in Section \ref{sec:double-dirichlet-series}. 

As the $\GL(3)$ $L$-function $L(s, \pi \otimes \chi_D)$ satisfies a functional equation $s \to 1-s$, the double Dirichlet series $\widetilde{Z}(s, w; \pi)$ is expected to satisfy an analogous functional equation. On the other hand, the Dirichlet series made of $\GL(1)$ $L$-functions obtained by  formal application of the quadratic reciprocity law, namely
\begin{equation}
\widehat{Z}(s, w; \pi) := \sum_{D, N} \frac{a_\pi(N)\chi_N(D)}{|N|^s|D|^w} =  \sum_{\substack{N}} \frac{L(w, \chi_N) a_\pi(N)}{|N|^s} ,
\end{equation}

\noindent is closely related to the original sum $\widetilde{Z}(s, w; \pi)$. Since the $\GL(1)$ $L$-function $L(w, \chi_N)$ satisfies a functional equation $w \to 1-w$, it is expected that $\widehat{Z}(s, w; \pi)$, and therefore $\widetilde{Z}(s, w; \pi)$, satisfies an analogous functional equation as well. However, this is not the case, as the quadratic characters only depend on the squarefree parts of the ideals $D$. As such, it is necessary to complete the double Dirichlet series via introducing correction factors $a(s, D, \pi)$ and $b(w, N, \pi)$ which are given by Dirichlet polynomials. The corrected double Dirichlet series is expected to take the form
\begin{equation}
Z(s, w; \pi) = \sum_{\substack{D}} \frac{L(s, \pi \otimes \chi_D)}{|D|^w} a(s, D, \pi) = \sum_N \frac{a_\pi(N) L(w, \chi_N)}{|N|^s} b(w, N, \pi)
\end{equation}

\noindent and to satisfy functional equations. These requirements yield necessary conditions on the desired correction factors, and Bump, Friedberg and Hoffstein \cite{bump_sums_2004} proved that these Dirichlet polynomials $a(s, D, \pi)$ and $b(w, N, \pi)$ exist and are unique in the case of $\GL(3)$. The two expected functional equations are then proven for the completed double Dirichlet series $Z(s, w; \pi)$ in Section \ref{sec:functiona$L$-equations}. They relate
\begin{align*}
Z(s, w; \pi)  \quad \text{and} \quad Z(\phi(s, w); \pi),  \\
Z(s, w; \pi)  \quad \text{and} \quad Z(\psi(s, w); \pi) ,
\end{align*}

\noindent where 
\begin{align*}
\phi(s, w) & = \left(1-s, w+3s-\tfrac{3}{2}\right), \\
\psi(s, w) & = \left(s+w-\tfrac{1}{2}, 1-w\right).
\end{align*}

We meromorphically continue the corrected double Dirichlet series $Z(s, w; \pi)$ to $\C^2$ in Section \ref{sec:meromorphic-continuation} by employing the functional equations of the series. To pull back the analytic properties of this corrected double Dirichlet series $Z(s, w; \pi)$ to the original double Dirichlet series $\widetilde{Z}(s, w; \pi)$, we carry out a technical sieving process in Section~\ref{sec:sieving-process}, expressing $\widetilde{Z}(s, w; \pi)$ as a twisted sum of functions of the form $Z(s, w; \pi) $. However, this sieving expresses $\widetilde{Z}(s, w; \pi)$ as an infinite sum of functions of type $Z(s, w; \pi)$, and the convergence of this sum may be too weak to ensure the expected meromorphic continuation. To this end, relevant estimates in the different aspects appearing in the summation are established in Section \ref{subsec:vertica$L$-bounds}, guaranteeing that $\widetilde{Z}(s, w; \pi)$ admits a meromorphic continuation around the point $(\frac{1}{2}, 1)$. It is in this process that the difficulties faced by Chinta and Diaconu have been overcome.

The proof then reduces to a computation of the residues $\widetilde{Z}(s, w; \pi)$ at $(\frac{1}{2}, 1)$, carried out in Section \ref{sec:loca$L$-analysis}. Indeed, after twisting $\pi$ by the character $\chi_N$ for a certain ideal $N$, the residue is proven to be an injective function of the $N$-th Fourier coefficient $a_\pi(N)$ of $\pi$ for large enough $N$. Therefore, if we assume that for every character $\chi_D \in X$,
\begin{equation}
L\left(\tfrac{1}{2}, \pi \otimes \chi_D\right) = \kappa \cdot  L\left(\tfrac{1}{2}, \pi' \otimes \chi_D\right),
\end{equation}

\noindent for a certain nonzero constant $\kappa$, then summing over $D$ yields an equality between the associated double Dirichlet series $\widetilde{Z}(\frac{1}{2}, w, \pi)$ and $\widetilde{Z}(\frac{1}{2}, w, \pi')$, and therefore between their residues. In particular, we can deduce the equality of almost all their Fourier coefficients by the injectivity mentioned above, and conclude by multiplicity one.

\subsection{Acknowledgements}
\label{subsec:acknowledgements}

We thank Farrell Brumley and Adrian Diaconu for enlightening discussions and Sun Yat-Sen University for its warm environment. The first-named author is supported in part by NSFC (No.11901585).

\section{Background}
\label{sec:background}

In this section, we first recap the definition of quadratic symbols on number fields and then state the functional equations of $L$-functions in the $\GL(1)$ and $\GL(3)$ settings.

\subsection{Quadratic symbols on number fields}
\label{subsec:quadratic-symbols}

Over $\Q$, the primitive quadratic characters \cite{montgomery_multiplicative_2006} are given by the Legendre symbol $\chi_d(n) = \legendre{d}{n}$ for $d$ varying among quadratic discriminants. Given a number field $F$, these characters have been generalized to the quadratic symbol of the form $\chi_d(N) = \legendre{d}{N}$ when $d \in F$ and $N$ is a fractional ideal of $F$ outside of ramified primes, see for instance \cite{cassels_algebraic_1967}. It is necessary for our purposes to extend the definition of these quadratic characters to symbols of the form $\chi_\d(N) = \legendre{\d}{N}$ for both $\d$ and $N$ integral ideals of $F$, without assuming $\d$ principal. This construction has been carried out by Fisher and Friedberg \cite[Section 1]{fisher_double_2004} for function fields, and adapted to number fields by Chinta, Friedberg and Hoffstein \cite[Section 2]{chinta_asymptotics_2005}. We provide here an account of their construction, adding details that we find useful and not explicitly available in the existing literature.

\begin{lem}
\label{lem:loca$L$-class-number}
There exists a finite set $S$ of finite places, such that $\mathcal{O}^{S}$ has class number one.
\end{lem}

\proof Let $h$ be the class number of $\O$. If $I$ is a non-principal ideal of $\O$, then $I^h$ is principal, so that $I^h = (a)$ for an $a \in F$. Let $S=\{\p_1,\ldots,\p_k\}$ be the finite set of prime ideals containing $(a)$, therefore containing $I$. In the Dedekind domain $\mathcal{O}^{S}$ of $S$-integers, the ideal $(a)$, and therefore also $I$, becomes trivial. The class group of $\mathcal{O}$ surjects onto that of $\mathcal{O}^{S}$, and the surjection kills
$I$, one of the nontrivial elements of the class group. Thus we get a Dedekind domain $\mathcal{O}^S$ with a smaller class number than $\mathcal{O}$. We finish by induction, adding at each step a finite number of primes to the set $S$, until the class number reaches one. \qed

Let $S$ be a finite set of places, containing the archimedean ones and the primes lying above 2. By Lemma \ref{lem:loca$L$-class-number}, it can be chosen such that $\mathcal{O}^{S_f}$ has class number one, where $S_f$ denotes the set of finite places in $S$. For an element $a$ of $F$, the Kronecker symbol $\legendre{a}{\cdot}$ is well-defined \cite{cassels_algebraic_1967} as the totally multiplicative function trivial on units and satisfying, for $\p$ a prime ideal of $F$,
\begin{equation}
\legendre{a}{\p} = a^{\frac{|\p| - 1}{\raisebox{-.3ex}{\tiny 2}}} \mod \ \p.
\end{equation}

It remains to extend the Kronecker symbol as a function of ideals in the numerator. For a place $v$ of $F$, let $F_v$ be the local completion at $v$. If $v$ is a finite place, let $\p_v$ denote the maximal ideal of $F_v$. Let $\c \in \O$ be defined by 
\begin{equation}
\label{def:c}
\c = \prod_{v \in S_f} \p_v^{n_v},
\end{equation}

where the exponents $n_v$ are chosen so that
\begin{equation}
 \left\{
\begin{array}{cl}
n_v = 1 & \text{if $v$ is not lying above $2$;}  \\
n_v & \text{such that } \mathrm{ord}_v(a-1) \geqslant n_v \Longrightarrow a \in F_v^{\times 2}.
\end{array}
 \right.
\end{equation}

Let $I(S)$ be the group of fractional ideals of $F$ prime to $S$ and $I^+(S)$ be the integral ideals contained in $I(S)$. Let $P_\c$  be the subgroup of $I(S)$ consisting of totally positive principal ideals $(a)$ with $a \equiv 1$ modulo $C$. Let $R_\c = I(S)/P_\c$ be the ray class group modulo $\c$, which is known to have finite cardinality. Let $H_\c = R_\c \otimes (\Z/2\Z) \simeq I(S) / I(S)^{2} P_\c$, so that $\widehat{H}_C$ consists of the quadratic characters of $R_C$. Since $H_\c$ is a finite abelian group, it decomposes as a direct sum of cyclic ones: there are elements $\bar{E}_1, \ldots, \bar{E}_k$ in $H_C$ such that
\begin{equation}
H_\c = \langle \overline{E}_1 \rangle \oplus \cdots \oplus \langle \overline{E}_k \rangle .
\end{equation}

Let $\mathcal{E}_0 = \{E_1, \ldots, E_k\}$ be a set of representatives for these generators, with $E_j \in I(S)$. Since $\mathcal{O}^{S_f}$ has class number one, \textit{i.e.} is a principal domain, for every $E_0 \in \mathcal{E}_0$ there is a generator $m_{E_0} \in \mathcal{O}^{S_f}$ such that $E_0 \O^{S_f} = m_{E_0} \mathcal{O}^{S_f}$. Let $\mathcal{E}$ be a full set of representatives for $R_\c$, namely the monoid generated by $\mathcal{E}_0$. More precisely, $\mathcal{E}$ is composed of the elements of the form
\begin{equation}
\label{def:decomposition-e}
E = \prod_{E_0 \in \mathcal{E}_0} E_0^{n_{E_0}} \quad \text{with} \quad n_{E_0} \geqslant 0. 
\end{equation}

For such an element $E$, define $m_E := \prod_{\mathcal{E}_0} m_{E_0}^{n_{E_0}}$, so that for every $E \in \mathcal{E}$, we get $E\mathcal{O}^{S_f} = m_E\mathcal{O}^{S_f}$. We can now extend the definition of the quadratic symbol. Any ideal $D$ in $I(S)$ decomposes as $(m)EG^2$ with $m \in F^\times$ such that $m \equiv 1 \mod \ \c$ (corresponding to the component in $P_\c$), $E \in \mathcal{E}$ a representative in $I(S)/ I(S)^2 P_\c$, and $G \in I(S)$ being the squareful part. We can now define the quadratic symbol appealing to the already known definition of $\legendre{\cdot}{D'}$ for elements of $F$ and $D' \in I(S)$ by letting, for all ideal $D \in I(S)$ decomposed as above,
\begin{equation}
\label{quadratic-symbo$L$-def}
\chi_D(D') = \legendre{D}{D'} := \legendre{m \cdot m_E}{D'}.
\end{equation}

All the good properties expected from a quadratic character are satisfied, as stated in the following proposition, coming from the very definition \eqref{quadratic-symbo$L$-def} and the properties of the classical quadratic symbol. It is analogous to the corresponding properties of Fisher and Friedberg \cite[Lemma 1.1]{fisher_double_2004} for function fields. 

\begin{prop}
\label{prop:properties-quadratic-symbol}
For $\d$ and $\d'$ integral ideals in $I^+(S)$, we have the properties:
\begin{itemize}
\item[(i)] $\chi$ is multiplicative: for every ideals $\d$ and $\d'$, we have $\chi_{\d\d'} = \chi_{\d}\chi_{\d'}$;
\item[(ii)] $\chi$ only depends on squarefree parts: if $\d$ has squarefree part $\d_0$, $\chi_\d = \chi_{\d_0}$;
\item[(iii)] Reciprocity law: for every $\d$ and $D'$ relatively prime ideals, $\chi_\d(\d')\chi_{\d'}(\d)$ is either $-1$ or $1$, depending only on the classes of $\d$ and $\d'$ in $H_C$. 
\end{itemize}
\end{prop}

The quadratic reciprocity law can be expressed as follows. For two relatively prime ideals $D$ and $D'$, define
\begin{equation}
\label{quadratic-reciprocity-law}
\eta(D, D') := \chi_D(D') \chi_{D'}(D).
\end{equation}

\noindent The value of $\eta(D, D')$ depends only on the classes of $D$ and $D'$ in $H_C$.  In particular we get
\begin{equation}
\chi_D(D') = \eta(D, D')\chi_{D'}(D).
\end{equation}

\emph{For the rest of this article, for $D \in I^+(S)$, we use $D_0$ to denote its square-free part and $D_1$ to denote its squareful part, i.e. $D = D_0 D_1^2$.}

\subsection{$L$-functions}
\label{subsec:$L$-functions}

\subsubsection{$\GL(1)$ case}

For $\d$ a squarefree ideal in $I^+(S)$, let $\chi_\d$ be the primitive ray class quadratic character over $F$ as defined in Section \ref{subsec:quadratic-symbols}. For every prime ideal $P$ and $w \in \C$, define the local $L$-factor by
\begin{equation}
L_P(w, \chi_D) := \left( 1 - \frac{\chi_D(P)}{|P|^{w}} \right)^{-1}.
\end{equation}

Let $L_f(w, \chi_\d)$ be the finite part of the associated $L$-function, defined by
\begin{equation}
L_f(w, \chi_D) := \prod_{P} L_P(w, \chi_D) = \sum_N \frac{\chi_D(N)}{|N|^w},
\end{equation}

where the product is over prime ideals $P$ of $\mathcal{O}$ and the sum over non zero ideals $N$ of $\O$. It converges on a right half-plane. The completed $L$-function 
\begin{equation}
\label{$L$-function-dirichlet-characters-completed}
L(w, \chi_\d) := \left( \frac{2^{r_1}|D_F|}{(2\pi)^d} \right)^{w/2} \Gamma\left( \frac{w}{2} \right)^{r_1} \Gamma(w)^{r_2} L_f(w, \chi_\d)
\end{equation}

\noindent is known to have analytic continuation to the whole complex plane, except for a simple pole at $w=1$ if $\chi_D$ is the trivial conductor \cite{goldmakher_quadratic_2013}. Moreover, the completed $L$-function satisfies the functional equation
\begin{equation}
\label{functiona$L$-equation-dirichlet-characters-completed}
L(w, \chi_\d) = \varepsilon(w, \chi_D) L(1-w, \chi_\d), 
\end{equation} 

\noindent where the $\varepsilon$-factor is defined by  \cite[Equation (3.1)]{goldmakher_quadratic_2013}
\begin{equation}
\label{epsilon-factor-dirichlet-characters}
\varepsilon(w, \chi_D) := |\d_0|^{\frac{1}{2}-w}.
\end{equation}

\noindent For all finite set $S$ of places, introduce the prime-to-$S$ partial $L$-function
\begin{equation}
L^S(w, \chi_D) := \prod_{\substack{P \in I^+(S) \\ P \text{ prime}}} L_P(w, \chi_D) = \sum_{N \in I^+(S)} \frac{\chi_D(N)}{|N|^w}.
\end{equation}

The functional equation \eqref{functiona$L$-equation-dirichlet-characters-completed} yields a functional equation for the partial $L$-function $L^S(w, \chi_\d)$:
\begin{equation}
\label{functiona$L$-equation-dirichlet-characters}
L^S(w, \chi_\d) = \varepsilon(w, \chi_D) L^S(1-w, \chi_\d) \prod_{\nu \in S} \frac{L_\nu(1-w, \chi_\d)}{L_\nu(w, \chi_\d)},
\end{equation} 

\noindent in which the last product above is a meromorphic function.

\subsubsection{$\GL(3)$ case}

Let $\pi$ be a self-contragredient cuspidal automorphic representation of $\GL(3)$ over $F$. Gelfand, Graev and Pyatetskii-Shapiro \cite{gelfand_representation_1969} introduced the finite-part $L$-function attached to $\pi$. For every prime ideal $P$, introduce the local $L$-factor
\begin{equation}
L_P(s, \pi) := \prod_{j=1}^3 \left( 1 - \frac{\gamma_j(\p)}{|\p|^{s}} \right)^{-1},
\end{equation}

\noindent where the $\gamma_j(\p)$ are the local spectral parameters of $\pi$, also called its Satake parameters when $\pi_\p$ is unramified. The best known bound towards the Ramanujan conjecture for $\GL(3)$ over general number fields is due to Blomer and Brumley \cite[Theorem 1]{blomer_ramanujan_2011} and is given by
\begin{equation}
\gamma_j(\p) \ll |\p|^{5/14 + \varepsilon}.
\end{equation}

Introduce for later purposes the convenient notation, for $M \in I^+(S)$ and $j \in \{1, 2, 3\}$,
\begin{equation}
\gamma_j(M) := \prod_{P^r || M} \gamma_j \left( P^r \right).
\end{equation}

\emph{Moreover, we make the convention that if the indices $j$ appear without being explicitly defined in a sum (resp. product), then it is understood that the expression is a sum (resp. product) over $j \in \{1, 2, 3\}$.}  The finite-part $L$-function attached to $\pi$ is defined \cite{cogdell_analytic_2004} by the Euler product of  local $L$-factors, namely
\begin{equation}
\label{$L$-function-representations-finite-part}
L_f(s, \pi) := \prod_P L_P(s, \pi) = \sum_{N} \frac{a_\pi(N)}{|N|^s},
\end{equation}

\noindent where the sum runs over nonzero integral ideals of $F$ and the product over prime ideals $P$. The $a_\pi(\n)$ are called the Fourier coefficients of $\pi$. For $\d \in I^+(S)$, the quadratic twist by $\chi_\d$ of its $L$-function is defined by
\begin{equation}
\label{$L$-function-twisted-representations-finite-part}
L_f(s, \pi \otimes \chi_\d) := \sum_{N} \frac{a_\pi(N)\chi_\d(N)}{|N|^s} = \prod_\p \prod_{j=1}^3 \left( 1 - \frac{\gamma_j(\p) \chi_\d(\p)}{|\p|^{s}} \right)^{-1}.
\end{equation}

Let $c(\pi)$ the arithmetic conductor of $\pi$, and let $\varepsilon_\pi$ be the root number of $\pi$ which is a complex number of modulus one depending only on $\pi$. The $L$-function above admits a completion $L(s, \pi \otimes \chi_\d)$, adding archimedean factors $L_v(s, \pi \otimes \chi_D)$ for $v | \infty$, made of explicit Euler gamma functions. This completed $L$-function is entire for cuspidal automorphic representation $\pi$ and satisfies the functional equation \cite[Equation (3.3)]{fisher_sums_2003} 
\begin{equation}
\label{functiona$L$-equation-twisted-representation-completed}
L(s, \pi \otimes \chi_\d) =\varepsilon(s, \pi \otimes \chi_D) L(1-s, \pi \otimes \chi_\d),
\end{equation}

\noindent where the $\varepsilon$-factor is defined by
\begin{equation}
\label{epsilon-factor-representation}
\varepsilon(s, \pi \otimes \chi_D)  :=  \varepsilon_\pi |\d_0|^{3\left(\frac{1}{2} - s\right)} c(\pi)^{ \frac{1}{2} - s}.
\end{equation}

\noindent For a finite set $S$ of places, the functional equation of the corresponding partial $L$-function is given by
\begin{align}
\label{functiona$L$-equation-twisted-representation}
 L^S(s, \pi \otimes \chi_\d) =  \varepsilon(s, \pi \otimes \chi_D) L^S(1-s, \pi \otimes\chi_\d) \prod_{\nu \in S} \frac{L_\nu(1-s, \pi \otimes \chi_D)}{L_\nu(s, \pi \otimes \chi_D)},
\end{align}

\noindent where the last product is a meromorphic function.

\section{Double Dirichlet series}
\label{sec:double-dirichlet-series}

\subsection{Pure and corrected double Dirichlet series}
\label{subsec:dds-pure-and-corrected}

Let $\pi$ be an automorphic cuspidal representation of $\GL(3)$ over $F$. Let $S$ be a finite set of places as chosen in Section \ref{subsec:quadratic-symbols}, so that the quadratic symbols $\chi_D$ are defined for $D \in I(S)$. From now on, all the summation variables appearing are assumed to be in $I^+(S)$, unless otherwise stated. Define the double Dirichlet series associated to the partial $L$-function $L^S(s, \pi \otimes \chi_D)$
\begin{equation}
\label{dds-pure-definition}
\widetilde{Z}^{S}(s, w; \pi) := \sum_{D} \frac{L^S(s, \pi \otimes \chi_D)}{|D|^w},
\end{equation}

which converges for $\sigma = \Re(s) > 1$ and  $\tau = \Re(w) > 1$. All the double Dirichlet series appearing in this article are convergent in this domain, and $s$ and $w$ lie in it unless otherwise stated. Writing the $L$-function attached to $\pi \otimes \chi_D$ as a Dirichlet series \eqref{$L$-function-twisted-representations-finite-part} yields the more explicit expression
\begin{equation}
\label{dds-expanded-definition}
\widetilde{Z}^{S}(s, w; \pi) = \sum_{D, M} \frac{a_\pi(M)\chi_D(M)}{|D|^w |M|^s}.
\end{equation}

Following \cite[Theorem 2.1]{bump_sums_2004}, for all characters $\alpha, \beta \in \widehat{H}_C$, there exist unique Dirichlet polynomials $a^S(s, D, \pi, \alpha)$ and $b^S(w, M, \pi, \beta)$ satisfying the conditions
\begin{align}
\sum_{D \in I^+(S)} \frac{L^S(s, \pi \otimes \chi_D \alpha)}{|D|^w}  a^S(s, D, \pi, \alpha) \beta(D) &= \sum_{M \in I^+(S)} \frac{L^S(w, \chi_M \beta)}{|M|^s }b^S(w, M, \pi, \beta) \alpha(M), \label{basic-identity}\\
\label{corrrection-factor-functiona$L$-equation}
 a^S(s, D, \pi, \alpha)  &= |D_1|^{3(1-2s)} \chi_\pi(D_1^2) a^S(1-s, D, {\pi}, \alpha), \\
\label{correction-factor-functiona$L$-equation-b}
 b^S(w, M, \pi, \beta)  &= |M_1|^{1-2s}  b^S(1-w, M, \pi, \beta).
\end{align}

Introduce $S_l = S\cup \{l\}$. Moreover, the coefficients also satisfy, see \cite[Equation (4.9)]{diaconu_multiple_2003}, 
\begin{equation}
\label{extraprop}
a^S(s, \pi, Dl, \alpha) = a^{S_l}(s, \pi , D, \alpha \chi_{l}).
\end{equation}

In particular, the relation \eqref{basic-identity} is the "basic identity" of Bump, Friedberg and Hoffstein. Moreover, from the explicit description \cite[Equations (2.1)-(2.3) pages 29-30 and 37]{bump_sums_2004} of these correction factors, we have the bounds for $\sigma, \tau > 1$,
\begin{align}
\label{correction-factor-bound-a}
a^S(s, D, \pi, \alpha) &  \ll_\varepsilon 1, \\
\label{correction-factor-bound-a-12}
a^S(1/2, D, \pi, \alpha) &  \ll_\varepsilon |D_1|^{5/7+\varepsilon}, \\
\label{correction-factor-bound-b}
b^S(w, M, \pi, \beta) &  \ll_\varepsilon a_\pi(M).
\end{align}

For convenience, we name the twisted corrected double Dirichlet series in \eqref{basic-identity}
\begin{align}
\label{dds-corrected-s-definition}
Z^{S}(s, w; \pi, \alpha, \beta) :=  \sum_{D} \frac{L^S(s, \pi \otimes \chi_D \alpha)}{|D|^w}a^S(s, D, \pi, \alpha) \beta(D)
\end{align}

\textit{Remark.} The correction factor $a^S(s, D, \pi, \alpha)$ is trivial for squarefree ideals $D$. As such, this corrected double Dirichlet series amounts to mollifying the terms for non-squarefree ideals in \eqref{dds-pure-definition}. The double Dirichlet series $Z^{S}(s, w; \pi, \alpha, \beta)$ and $\widetilde{Z}^{S}(s, w; \pi)$ are therefore expected to be closely related: this is indeed the case, as shown in Section \ref{sec:sieving-process}.

\subsection{Ray class selection}
\label{subsec:quadratic-reciprocity}

Let $E$ be a representative of a ray class in $H_C$, and $\delta_E$ the characteristic function of this class, in the sense that for every ideal $D$, $\delta_E(D)$ takes value $1$ if $D$ is in the same class as $E$ modulo $I(S)^2 P_C$, and $0$ otherwise. The partition of $H_C$ in classes yields
\begin{equation}
1 = \sum_{[E] \in H_C} \delta_E,
\end{equation} 

\noindent where the summation runs over a set of representatives of classes in $H_C$. Moreover, letting $h_C = |H_C|$, the orthogonality relations for characters in $H_C$ can be rephrased as follows: for every element $E$ in $H_C$, $\alpha,\beta$ in $\widehat{H}_C$,
\begin{align}
\label{orthogonality-relation}
h_C^{-1} \sum_{\rho \in \widehat{H}_C} \rho(E)^{-1} \rho &= \delta_E, \\
\label{partition-classes}
h_C^{-1} \sum_{[E'] \in H_C} \alpha(E')^{-1} \beta(E') &= \begin{cases} 1, &\text{ if $\alpha= \beta$,} \\ 0, &\text{ otherwise}. \end{cases}
\end{align}

To allow for $\alpha,\beta$ to be taken as functions on $H_C$ in the double Dirichlet series, we extend the definition of the double Dirichlet series. To be precise, for any functions $\delta, \delta'$ on $H_C$ written as sums of characters
\begin{align*}
\delta & = \sum_{\rho \in \widehat{H}_C} \lambda_\rho \cdot \rho, \\
\delta' & = \sum_{\rho' \in \widehat{H}_C} \mu_{\rho'} \cdot \rho',
\end{align*}

define

\begin{align*}
Z^S(s, w; \pi, \delta , \delta') &: = \sum_{\rho, \rho' \in \widehat{H}_C} \lambda_\rho \mu_{\rho'} \cdot  Z^S(s, w; \pi, \rho , \rho'),
\end{align*}

Under the extension, relations \eqref{orthogonality-relation} and \eqref{partition-classes} translate into the following:

\begin{lem}[Sieving by classes in $H_C$]
\label{lem:sieving-classes-Hc}
Let $E, E'$ be ideals of $I(S)$ and $\alpha,\beta \in \widehat{H}_C$. We have
\begin{align*}
Z^S(s, w; \pi, \alpha \delta_E, \beta \delta_{E'}) & = h_C^{-2} \sum_{\rho, \rho' \in \widehat{H}_C} \rho(E)^{-1}\rho'(E')^{-1} \cdot  Z^S(s, w; \pi, \alpha \rho, \beta \rho'), \\[0.5em]
Z^S(s, w; \pi, \alpha, \beta) & = \sum_{[E], [E'] \in H_C} Z^S(s, w; \pi, \alpha \delta_E, \beta \delta_{E'}),
\end{align*}

where this last sum is over a set of representatives in $H_C$. 
\end{lem}

This sieving lemma allows us to select a single class for any summation variable appearing in the double Dirichlet series at the expense of a finite linear combination. Separating into different ray classes is particularly useful when we apply reciprocity laws on the quadratic character. Also note that the linear combination appearing in the above lemma has length and coefficients uniformly bounded in $\pi$, $D$, $M$, $\alpha$ and $\beta$, and therefore enjoy the same convergence properties as the full double Dirichlet series. Henceforth, linear combinations of this sort will be referred to as \emph{uniformly finite linear combinations}.

\section{Functional equations}
\label{sec:functiona$L$-equations}

\subsection{Epsilon-factors}
\label{subsec:functiona$L$-equation-a}

Let $\pi$ be a self-contragredient cuspidal automorphic representation of $\GL(3, \A)$ with central character $\chi_\pi$, unramified outside $S$ and moreover assumed to be a principal series representation at every place outside $S$. It is possible to suppose so without loss of generality \cite{bump_automorphic_1997}, since this is the case except for a finite number of places that can be added to $S$. 
\begin{lem}[Class expression for $\varepsilon$-factors]
\label{lem:epsilon-factor-relation-DE}
Let $D, E$ be integral ideals prime to $S$. Suppose that $D$ and $E$ lie in the same class in $H_C$, implying $\chi_D = \chi_E \chi_m$ for a certain $m \in F^\times$ satisfying $m \equiv 1 \ (\mathrm{mod} \ C)$.  Let $\pi$ be a self-contragredient cuspidal automorphic representation of $\GL(n, \A)$ for $n=1$ or $n=3$, unramified outside $S$. Then the value of $\varepsilon_v(s, \pi_v \otimes \chi_{D, v}) $ depends solely on the class of $D$ in $H_C$ for places $v \in S$. Moreover, 
\begin{equation}
\label{epsilon-factor-relation-DE}
\varepsilon(s, \pi \otimes \chi_D) = \chi_\pi\left(\frac{D_0}{E_0}\right) \left|\frac{D_0}{E_0}\right|^{n(\frac{1}{2}-s)} \varepsilon(s, \pi \otimes \chi_E).
\end{equation}
\end{lem}

\proof Note that the result is obvious if the central values are zero, so we can assume they are nonzero until the end of the proof. By the definition \eqref{epsilon-factor-representation} of the $\varepsilon$-factor, 
\begin{equation}
\varepsilon(s, \pi \otimes \chi_D) = \left| c(\pi \otimes \chi_D)\right|^{\frac{1}{2}-s}  \varepsilon\left(\tfrac{1}{2}, \pi \otimes \chi_D\right).
\end{equation}

Let $v \in S$. Recall that, thanks to the decomposition $D=(m)EG^2$ obtained in Section \ref{subsec:quadratic-symbols}, we have $\chi_D = \chi_E\chi_m$. It is assumed that $v_P(m-1) > 0$ for finite places $P \in S$, that is to say $m \equiv 1$ modulo $C$. Since $C$ only consists of finite primes in $S$, we deduce that $\chi_{D, v} = \chi_{E, v}$ for $v \in S$. 

Let $v \notin S$. In that case, the assumption that $\pi_v$ is an unramified principal series implies $\pi = \psi_1$ if $n=1$ or $\pi = \pi(\psi_1, \psi_2, \psi_3)$ if $n=3$, for characters $\psi_j$ of $F_v^\times$. The $\varepsilon$-factors \eqref{epsilon-factor-dirichlet-characters} decompose as a product of local $\varepsilon$-factors of the form. Let $\theta_v$ an additive character on $F_v$. For any $j \in \{1, 2, 3\}$,
\begin{equation}
\varepsilon_v\left( \tfrac{1}{2}, \psi_{j}\chi_{D, v}, \theta_v \right) = \psi_j(P_v)^{\mathrm{ord}_v(c(\chi_D))} \varepsilon_v\left( \tfrac{1}{2},\chi_{D, v}, \theta_v \right).
\end{equation}

The same formula holds for the specific choice $D=E$, so that in particular we get by dividing these nonzero quantities and taking product of $j$, and since $\psi_1\psi_2\psi_3 = \chi_\pi$,
\begin{equation}
\frac{\varepsilon\left(\frac{1}{2}, \pi \otimes \chi_D\right)}{\varepsilon\left(\frac{1}{2}, \pi \otimes \chi_E\right)} = \chi_\pi \left( \frac{D_0}{E_0} \right)  \frac{\varepsilon\left(\frac{1}{2}, \chi_D\right)^n}{\varepsilon\left(\frac{1}{2},  \chi_E\right)^n} .
\end{equation}

The central values for quadratic characters are assumed to be nonzero, so the ratio of the epsilon factors on the right is 1. With the fact that
\begin{equation}
\frac{c(\pi \otimes \chi_D)}{c(\pi \otimes \chi_E)} = \left|\frac{D_0}{E_0}\right|^n,
\end{equation}

\noindent the result is established. \qed

This lemma explicitly gives the dependence on a fixed representative in $H_C$.  We immediately deduce from this and the fact that $(\pi \otimes \chi_M) \otimes \chi_D = \pi \otimes \chi_{DM}$ the following corollary, which is \cite[Corollary 1.11]{li_determination_2007}.

\begin{coro}
\label{coro:varepsilon-twist-by-character-relation-DE}
Let $L \in I^+(S)$ squarefree. Suppose $D, E \in I^+(S)$ are in the same class of $H_C$. Let $\pi$ be a self-contragradient cuspidal automorphic representation of $\GL(n, \A)$ for $n=1$ or $n=3$, unramified outside $S$. Then
\begin{equation}
\varepsilon(s, (\pi \otimes \chi_L)\otimes \chi_D) = \chi_\pi\left(\frac{D_0}{E_0}\right) \left| \frac{D_0}{E_0} \right|^{n\left(\frac{1}{2}-s\right)} \varepsilon(s, (\pi \otimes \chi_L)\otimes \chi_E).
\end{equation}
\end{coro}

\subsection{Functional equation in $s$}
\label{subsec:functiona$L$-equation-s}

As the $L$-function $L(s, \pi)$ satisfies a functional equation \eqref{functiona$L$-equation-twisted-representation-completed} relating $s \leftrightarrow 1-s$, the associated double Dirichlet series is expected to satisfy a similar functional equation.	Recall that 
\begin{equation}
\phi(s, w) = \left(1-s, w+3s-\tfrac{3}{2}\right).
\end{equation}

\begin{prop}[Functional equation $\phi  :s \leftrightarrow 1-s$]
\label{prop:fe-s}
For a self-contragradient cuspidal automorphic representation $\pi$ on $\GL(3)$ with central character $\chi_\pi$, $\alpha, \beta \in \widehat{H}_C$, $E$ a class in $H_C$, the double Dirichlet series $Z^S(s, w; \pi, \alpha, \beta \delta_E)$ satisfies the functional equation for $\sigma, \tau > 1$,
\begin{equation}
\label{fe1}
 Z^S(s, w; \pi, \alpha,  \beta \delta_E)  = Z^S\left(\phi(s, w) ; \pi, \alpha, \beta \chi_\pi \delta_E \right) \frac{\varepsilon(s, \pi \otimes \chi_E \alpha)}{\chi_\pi(E_0) |E_0|^{3(\frac{1}{2}-s)}} \prod_{\nu \in S} \frac{L_\nu(1-s, \pi \otimes \chi_E \alpha)}{L_\nu(s, \pi \otimes \chi_E \alpha)}.
\end{equation}
\end{prop}

\proof This is essentially summing over $D$ the functional equations of the involved $L$-functions, $\varepsilon$-factors  and correction factors. We follow the steps of Fisher and Friedberg \cite[Theorem 3.1]{fisher_double_2004}. By \eqref{corrrection-factor-functiona$L$-equation}, the Dirichlet polynomial satisfies the functional equation
\begin{equation}
\label{fe-a}
a^S(s, D, \pi, \alpha) = a^S(1-s, D, \pi, \alpha) \chi_{\pi}(D_1^2) |D_1|^{6\left(\frac{1}{2}-s\right)}.
\end{equation}

By Lemma \ref{lem:epsilon-factor-relation-DE}, the $\varepsilon$-factors at places in $S$ only depend on the class $E$ of $D$, not on $D$ itself. Hence, we deduce that
\begin{align}
 L^S(s, \pi \otimes \chi_D \alpha) &=  L^S(1-s, \pi \otimes \chi_D \alpha)  
\varepsilon(s, \pi \otimes \chi_D \alpha) \prod_{v \in S} \frac{L_v(1-s, \pi \otimes \chi_D \alpha)}{L_v(s, \pi \otimes \chi_D \alpha)} \notag \\
&= L^S(1-s, \pi \otimes \chi_D \alpha)  
\varepsilon(s, \pi \otimes \chi_D \alpha) \prod_{v \in S} \frac{L_v(1-s, \pi \otimes \chi_E \alpha)}{L_v(s, \pi \otimes \chi_E \alpha)}.
\end{align}

Multiplying the above with the Dirichlet polynomial, we obtain
\begin{equation}
\label{aaa}
\begin{split}
L^S(s, \pi \otimes \chi_D \alpha) a^S(s, D, \pi, \alpha) &=  L^S(1-s, \pi \otimes \chi_D \alpha)  a^S(1-s, D, \pi, \alpha) \\
&\qquad \times |D_1|^{6\left(\frac{1}{2}-s\right)} \chi_\pi(D_1^2) \varepsilon(s, \pi \otimes \chi_D \alpha) \prod_{v \in S} \frac{L_v(1-s, \pi \otimes \chi_E \alpha)}{L_v(s, \pi \otimes \chi_E \alpha)}.
\end{split}
\end{equation}

Furthermore, incorporating the expression of the $\varepsilon$-factor in \eqref{epsilon-factor-relation-DE}, we have
\begin{equation}
\label{blah}
\begin{split}
L^S(s, \pi \otimes \chi_D \alpha) a^S(s, D, \pi, \alpha) &= L^S(1-s, \pi\otimes \chi_D \alpha) a^S(1-s, D, \pi, \alpha)  \\
&  \qquad \times \frac{\varepsilon(s, \pi \otimes\chi_E \alpha)}{\chi_\pi(E_0)|E_0|^{3(\frac{1}{2}-s)}} \chi_\pi(D) |D|^{3\left(\frac{1}{2}-s\right)} \prod_{v \in S} \frac{L_v(1-s, \pi \otimes \chi_E \alpha)}{L_v(s, \pi \otimes \chi_E \alpha)}.
\end{split}
\end{equation}

In order to recover the double Dirichlet series $Z^S(s, w; \pi,\alpha,  \beta \delta_E)$, we need to sum over $D$. Due to the presence of $\delta_E$ in the double Dirichlet series, we sum over all the ideals in the $H_C$-class $[E]$.  After multiplying by $\beta(D) |D|^{-w}$ and summing \eqref{blah} over $D \in [E]$, for $w$ with large enough real part, we get
\begin{align*}
Z^S(s, w ; \pi, \alpha, \beta \delta_E) &= \sum_{D \in [E]} \frac{L^S(s, \pi \otimes \chi_D \alpha)}{|D|^w}  a^S(s, D, \pi, \alpha) \beta(D) \\
&= \frac{\varepsilon(s, \pi \otimes\chi_E \alpha)}{\chi_\pi(E_0) |E_0|^{3(\frac{1}{2}-s)}}  \prod_{v \in S} \frac{L_v(1-s, \pi \otimes \chi_E \alpha)}{L_v(s, \pi \otimes \chi_E \alpha)}   \\
& \qquad \qquad \times  \sum_{D \in [E]}  \frac{L^S(1-s, \pi\otimes \chi_D \alpha) }{|D|^{w + 3s - 3/2}} a^S(1-s, D, \pi, \alpha) \beta(D) \chi_\pi(D).
\end{align*} 
 
Recognizing the series in the last equality as $Z^S(\phi(s, w); \pi, \alpha, \beta \chi_\pi \delta_E)$ finishes the proof. \qed

Summing over the classes $E$ in $H_C$ therefore yields an explicit functional equation relating $Z^S(s, w; \pi, \alpha, \beta)$ and $Z^S(\phi(s, w); \pi, \alpha, \beta)$. Since these sums are uniformly finite by Lemma \ref{lem:sieving-classes-Hc}, the analytic behavior is similar and we can afford to switch freely between statements on $Z^S(s, w; \pi, \alpha \delta_E, \beta \delta_{E'})$ and statements on $Z^S(s, w; \pi, \alpha, \beta)$. 
%


For later sections, it is necessary to make explicit the dependence when further skipping primes of a chosen ideal $r$. For this purpose, we have the following refined version of the functional equation.

\begin{prop}[Refined functional equation $\phi : s \leftrightarrow 1-s$]
\label{prop:fe-s-refined}
For a self-contragradient cuspidal automorphic representation $\pi$ on $\GL(3)$ with central character $\chi_\pi$, $\alpha, \beta \in \widehat{H}_C$, $E$ a class in $H_C$ out of $S$, and a squarefree ideal $r \in I^+(S)$. Let $S_r$ be the finite set of places in $S$ along with those dividing $r$, and $f_r(\alpha)$ be the product of primes dividing $r$ and the conductor of $\alpha$. We view $\delta_E$ as a function on $I^+(S_r)$. We have that
\begin{equation}
\prod_{P | r /f_r(\alpha)} \prod_j \left(1- \frac{ \alpha \gamma_j(P)^2}{|P|^{2-2s}}  \right)  Z^{S_r}(s, w; \pi, \alpha, \beta \delta_E)
\end{equation}

is a uniformly finite linear combination of expressions of the type, for $\rho, \rho' \in \widehat{H}_C$,
\begin{align*}
& \sum_{l_j  | r/f_r(\alpha)} \frac{ \alpha \gamma_j(l_j)}{|l_j|^{1-s}}  \sum_{m_j  | r /f_r(\alpha)}  \frac{\mu \alpha \gamma_j(m_j)}{|m_j|^{s}} \rho(l_1l_2l_3m_1m_2m_3) \prod_{\nu \in S} \frac{L_\nu(1-s, \pi \otimes \chi_E \alpha)}{L_\nu(s, \pi \otimes \chi_E \alpha)}  \\
& \qquad \times   \frac{\varepsilon(s, \pi \otimes \chi_E \alpha)}{\chi_\pi(E_0)|E_0|^{3(\frac{1}{2}-s)}} Z^{S_r}(\phi(s, w) ; \pi, \alpha, \beta \rho' \chi_\pi \chi_{l_1l_2l_3m_1m_2m_3}).
\end{align*}
\end{prop}

\proof The refinement compared to the functional equation obtained in Proposition \ref{prop:fe-s} consists of undisclosing the $L$-factors for the places $P$ dividing $r$, according to \eqref{functiona$L$-equation-twisted-representation}. Similar to the previous proof, we see that
\begin{align*}
  Z^{S_r}(s,w;\pi,\alpha,\beta \delta_E) = &\frac{\varepsilon(s,\pi \otimes \chi_E \alpha)}{\chi_\pi(E_0) |E_0|^{3(1/2-s)}}  \prod_{v \in S} \frac{L_v(1-s, \pi \otimes \chi_E \alpha)}{L_v(s, \pi \otimes \chi_E \alpha)}  \\
  &\times \sum_{D \in [E]}  \frac{L^{S_r}(1-s, \pi\otimes \chi_D \alpha) }{|D|^{w + 3s - 3/2}} a^{S_r}(1-s, D, \pi, \alpha) \beta \chi_\pi(D) \prod_{P | r} \frac{L_P(1-s, \pi \otimes \chi_D \alpha)}{L_P(s, \pi \otimes \chi_D \alpha)}
\end{align*}
The product over prime factors $P$ of $r$ can be changed to prime factors of $r/f_r(\alpha)$, due to the character $\alpha$ being applied to the prime ideal in the $L$-factor. Noting that $\chi_D(P)^2 = 1$, we see that 
\begin{equation}
  1- \frac{\alpha\gamma_j(P)^2}{ |P|^{2-2s}} = \left(1- \frac{\alpha \gamma_j \chi_D(P)}{ |P|^{1-s}} \right) \left(1+\frac{\alpha \gamma_j \chi_D(P)}{ |P|^{1-s}}\right).
\end{equation}
With this observation, we see that
\begin{align*}
  & \prod_{P | r /f_r(\alpha)} \prod_j \left(1- \frac{ \alpha \gamma_j(P)^2}{|P|^{2-2s}}  \right)  Z^{S_r}(s, w; \pi, \alpha, \beta \delta_E) \\
  =& \frac{\varepsilon(s,\pi \otimes \chi_E \alpha)}{\chi_\pi(E_0) |E_0|^{3(1/2-s)}}  \prod_{v \in S} \frac{L_v(1-s, \pi \otimes \chi_E \alpha)}{L_v(s, \pi \otimes \chi_E \alpha)} \sum_{D \in [E]}  \frac{L^{S_r}(1-s, \pi\otimes \chi_D \alpha) }{|D|^{w + 3s - 3/2}} \\
  &\times a^{S_r}(1-s, D, \pi, \alpha) \beta \chi_\pi(D) \prod_{P | r/f_r(\alpha)} \prod_j  \left( 1 - \frac{ \alpha \gamma_j \chi_D(P)}{|P|^s} \right) \left( 1 + \frac{\alpha \gamma_j \chi_D(P)}{ |P|^{1-s}} \right) \\
 = & \sum_{l_j  | r/f_r(\alpha)} \frac{ \alpha \gamma_j(l_j)}{|l_j|^{1-s}}  \sum_{m_j  | r /f_r(\alpha)}  \frac{\mu \alpha \gamma_j(m_j)}{|m_j|^{s}} \frac{\varepsilon(s,\pi \otimes \chi_E \alpha)}{\chi_\pi(E_0) |E_0|^{3(1/2-s)}} \prod_{\nu \in S} \frac{L_\nu(1-s, \pi \otimes \chi_E \alpha)}{L_\nu(s, \pi \otimes \chi_E \alpha)} \\
 & \times \sum_{D} \frac{L^{S_r}(1-s, \pi\otimes \chi_D \alpha) }{|D|^{w + 3s - 3/2}} a^{S_r}(1-s, D, \pi, \alpha) \beta \chi_\pi \delta_E(D) \chi_D(l_1l_2l_3 m_1m_2m_3) 
\end{align*}
It remains to integrate the quadratic symbol  $\chi_D(l_1l_2l_3m_1m_2m_3)$ in the sum over $D$, seeing it as a character in $D$ in order to recognize a genuine double Dirichlet series. The quadratic reciprocity law for quadratic characters \eqref{quadratic-reciprocity-law} writes, for $D$ in the same class as $E$,
\begin{equation}
\chi_D(l_1l_2l_3m_1m_2m_3) = \eta(E, l_1l_2l_3m_1m_2m_3) \chi_{l_1l_2l_3m_1m_2m_3} (D).
\end{equation}

Thus, sieving the classes of $l_1, \cdots, m_3$ in $H_C$ by the orthogonality relations \eqref{orthogonality-relation}, we are reduced to a uniformly finite linear combination of sums of the form
\begin{align*}
& \sum_{l_j | r/ f_r(\alpha)} \frac{\alpha \gamma_j(l_j)}{|l_j|^{1-s}}  \sum_{m_j | r/ f_r(\alpha)}  \frac{\mu \alpha \gamma_j(m_j)}{|m_j|^{s}} \frac{\varepsilon(s, \pi \otimes \chi_E)}{\chi_\pi(E_0)^{-1}|E_0|^{3(\frac{1}{2}-s)}}  \prod_{\nu \in S} \frac{L_\nu(1-s, \pi \otimes \chi_E)}{L_\nu(s, \pi \otimes \chi_E)}  \\
& \qquad \times  \rho(l_1l_2l_3m_1m_2m_3) Z^{S_r} \left(1-s, w+3s-\tfrac{3}{2} ; \pi , \alpha , \beta \rho' \chi_\pi \chi_{l_1l_2l_3m_1m_2m_3}\right),
\end{align*}

\noindent  for $\rho, \rho' \in \widehat{H}_C$. This is exactly the claimed result. \qed

\subsection{Functional equation in $w$}
\label{subsec:functiona$L$-equation-w}

The $L$-function $L(w, \chi_D)$ satisfies a functional equation \eqref{functiona$L$-equation-dirichlet-characters} relating $w \leftrightarrow 1-w$, so that the same is expected for the associated double Dirichlet series. This is indeed the case by the following proposition. Recall that
\begin{equation}
\psi(s, w) = \left(s+w-\tfrac{1}{2}, 1-w\right).
\end{equation}

\begin{prop}[Functional equation $ \psi :w \leftrightarrow 1-w$]
\label{prop:fe-w}
Consider an automorphic cuspidal representation $\pi$ on $\GL(3)$, $\alpha, \beta \in \widehat{H}_C$, $E$ a class in $H_C$ out of $S$. The double Dirichlet series  $Z^S(s, w; \pi, \alpha \delta_E, \beta)$ satisfies the functional equation
\begin{align}
\label{fe-w}
 Z^S(s, w; \pi, \alpha \delta_E, \beta )  & = Z^S(\psi(s, w) ; \pi, \alpha \beta \delta_E, \beta) \frac{\varepsilon(w, \chi_E \beta)}{|E_0|^{\frac{1}{2}-w} \beta(E_0)} \prod_{\nu \in S} \frac{L_\nu(1-w, \chi_E \beta)}{L_\nu(w, \chi_E \beta)}.
\end{align}
\end{prop}

\proof The proof is analogous to the one of Proposition \ref{prop:fe-s} and amounts to summing the functional equations for the product of $L$-functions $L(w, \chi_M \beta)$ and the correction factors $b^S(w, M, \beta, \pi)$. By \eqref{basic-identity}, 
\begin{equation}
Z^{S}(s, w; \pi, \alpha \delta_E, \beta) = \sum_{M \in [E]} \frac{L^S(w, \chi_M \beta) }{|M|^s} \cdot b^S(w, M, \pi, \beta) \alpha(M).
\end{equation}

Introducing the functional equations \eqref{functiona$L$-equation-dirichlet-characters-completed} for the $L$-function $L(w, \chi_M)$ and \eqref{correction-factor-functiona$L$-equation-b} for $b^S(w, M, \pi, \beta)$, and then using Lemma \ref{lem:epsilon-factor-relation-DE} to formulate the $\varepsilon$-factor in function of the representative $E$, we get
\begin{align*}
& L^S(w, \chi_M \beta) b^S(w, M, \pi, \beta)  \\
  = &L^S(1-w, \chi_M \beta) b^S(1-w, M, \pi, \beta) \prod_{\nu \in S} \frac{L_\nu(1-w, \chi_M \beta)}{L_\nu(w, \chi_M \beta)} \varepsilon(w, \chi_M\beta) |M_1|^{1-2w} \beta(M_1^2) \\
 = &L^S(1-w, \chi_M \beta) b^S(1-w, M, \pi, \beta)  \varepsilon(w, \chi_E \beta)  \beta\left(\frac{M_0}{E_0}\right) \left|\frac{M_0}{E_0}\right|^{\frac{1}{2}-w} |M_1|^{1-2w} \beta(M_1^2)  \prod_{\nu \in S} \frac{L_\nu(1-w, \chi_E \beta)}{L_\nu(w, \chi_E \beta)} .
\end{align*}

Multiplying by $\alpha(M)|M|^{-s}$ and summing over $M$ yields 
\begin{align*}
Z^{S}(s, w; \pi, \alpha \delta_E, \beta) & = \sum_{M \in [E]} \frac{L^S(1-w, \chi_M\beta) }{|M|^{s+w-\frac{1}{2}}} b^S(1-w, M, \pi, \beta)    \prod_{\nu \in S} \frac{L_\nu(1-w, \chi_E \beta)}{L_\nu(w, \chi_E \beta)}  \frac{\varepsilon(w, \chi_E\beta)}{|E_0|^{\frac{1}{2}-w} \beta(E_0)} \beta \alpha(M).
\end{align*}

Therefore we recognize the double Dirichlet series
\begin{align*}
Z^{S}(s, w; \pi, \alpha \delta_E, \beta) & =   Z^{S}\left(s+w-\tfrac{1}{2}, 1-w; \pi, \alpha \beta \delta_E, \beta\right) \frac{\varepsilon(w, \chi_E\beta)}{|E_0|^{\frac{1}{2}-w} \beta(E_0)} \prod_{\nu \in S} \frac{L_\nu(1-w, \chi_E \beta)}{L_\nu(w, \chi_E \beta)}.
\end{align*} 

\noindent which is exactly the claimed result. \qed

Analogous to Proposition \ref{prop:fe-s-refined}, there is a refined version of this functional equation in order to have an explicit dependence in a fixed squarefree ideal $r \in I^+(S)$.
\begin{prop}[Refined functional equation $\psi : w \leftrightarrow 1-w$]
\label{prop:fe-w-refined}
Consider an automorphic cuspidal representation $\pi$ on $\GL(3)$, $\alpha, \beta \in \widehat{H}_C$, $E$ a class in $H_C$ out of $S$, and a squarefree ideal $r \in I^+(S)$. Let $S_r$ be the finite set of places in $S$ along with those dividing $r$, and $f_r(\beta)$ be the product of primes dividing $r$ and the conductor of $\beta$. We view $\delta_E$ as a function on $I^+(S_r)$. We have that
\begin{equation}
 \prod_{P|r/f_r(\beta)}  \left(1-\frac{\beta(P)^2}{|P|^{2-2w}}\right) Z^{S_r}(s, w; \pi, \alpha \delta_E,\beta)
\end{equation}

can be written as a uniformly finite linear combination of expressions of the form
\begin{align*}
  \sum_{l|r/f_r(\beta)} \frac{\mu \beta(l)}{|l|^{w}}  &\sum_{m | r/f_r(\beta)} \frac{\beta(m)}{|m|^{1-w} } \frac{\varepsilon(w, \chi_E\beta)}{|E_0|^{\frac{1}{2}-w} \beta(E_0)} \prod_{\nu \in S} \frac{L_\nu(1-w, \chi_E \beta)}{L_\nu(w, \chi_E \beta)} \rho(lm) Z^{S_r}\left(\psi(s,w); \pi, \alpha \beta  \chi_{lm}\rho', \beta\right),
\end{align*}

where $\rho, \rho' \in \widehat{H}_C$.
\end{prop}

\proof  The proof is analogous to the one of Proposition \ref{prop:fe-s-refined}, and essentially amounts to carefully summing the associated functional equations \eqref{fe1} over all $M$ and compensate appropriately at places $P|r$. We have
\begin{equation}
\label{up}
Z^{S_r}(s, w; \pi, \alpha \delta_E, \beta) = \sum_{M \in [E]} \frac{L^{S_r}(w, \chi_M \beta)}{|M|^s} b^{S_r}(w, M, \pi, \beta) \alpha(M).
\end{equation}

The functional equation \eqref{functiona$L$-equation-dirichlet-characters-completed} gives
\begin{align*}
& L^{S_r}(w, \chi_M \beta) =\varepsilon(w, \chi_M\beta) L^{S_r}(1-w, \chi_M \beta) \prod_{P | r} \frac{1-\chi_M\beta(P) |P|^{-w}}{1-\chi_M \beta(P) |P|^{-1+w}} \prod_{\nu \in S} \frac{L_\nu(1-w, \chi_M\beta)}{L_\nu(w, \chi_M\beta)}.
\end{align*}

Recall that by \eqref{correction-factor-functiona$L$-equation-b} the local factors $b^{S_r}(w, N, \pi, \beta)$ satisfy the function equation
\begin{equation}
b^{S_r}(w, M, \pi, \beta) = |M_1|^{1-2w} \beta(M_1^2) b^{S_r}(1-w, M, \pi, \beta).
\end{equation}

Introducing these two functional equations in the expression \eqref{up} leads to
\begin{align*}
& Z^{S_r}(s, w; \pi, \alpha \delta_E, \beta)   \\
= &\frac{\varepsilon(w, \chi_E\beta)}{|E_0|^{\frac{1}{2}-w} \beta(E_0)} \prod_{\nu \in S} \frac{L_\nu(1-w, \chi_E \beta)}{L_\nu(w, \chi_E \beta)} \sum_{M \in [E]} \frac{L^{S_r}(1-w, \chi_M\beta)}{|M|^{s+w-\frac{1}{2}}} b^{S_r}(1-w, M, \pi, \beta) \\
& \qquad \qquad \times \prod_{P | r} \left(1-\frac{\chi_M\beta(P)}{ |P|^{w}}\right) \left( 1- \frac{\chi_M \beta(P)}{ |P|^{1-w}} \right)^{-1} \alpha \beta(M).
\end{align*}

Analogous to Proposition \ref{prop:fe-s-refined}, multiplying both sides by $\prod (1-\beta(P)^2|P|^{-2+2w})$, where the product is over primes $P$ dividing $r/f_r(\beta)$, gives
\begin{align*}
& \prod_{P|r/f_r(\beta)} \left(1-\frac{\beta(P)^2}{|P|^{2-2w}} \right) Z^{S_r}(s, w; \pi, \alpha \delta_E, \beta)  \\
& = \frac{\varepsilon(w, \chi_E\beta)}{|E_0|^{\frac{1}{2}-w} \beta(E_0)} \prod_{\nu \in S} \frac{L_\nu(1-w, \chi_E \beta)}{L_\nu(w, \chi_E \beta)} \sum_{M \in [E]} \frac{L^{S_r}(1-w, \chi_M\beta)}{|M|^{s+w-\frac{1}{2}}} b^{S_r}(1-w, M, \pi, \beta) \\
& \qquad \qquad \times \prod_{P | r/f_r(\beta)} \left(1-\frac{\chi_M\beta(P)}{ |P|^{w}}\right) \left( 1+\frac{\chi_M \beta(P)}{ |P|^{1-w}} \right) \alpha \beta(M) \\
& = \sum_{l|r/f_r(\beta)} \frac{\mu \beta(l)}{|l|^{w}}  \sum_{m | r/f_r(\beta)} \frac{\beta(m)}{|m|^{1-w} } \frac{\varepsilon(w, \chi_E\beta)}{|E_0|^{\frac{1}{2}-w} \beta(E_0)} \prod_{\nu \in S} \frac{L_\nu(1-w, \chi_E \beta)}{L_\nu(w, \chi_E \beta)} \\
&\qquad \qquad  \times  \sum_{M \in [E]} \frac{L^{S_r}(1-w, \chi_M\beta)}{|M|^{s+w-\frac{1}{2}}} b^{S_r}(1-w, M, \pi, \beta) \alpha \beta(M) \chi_M(lm)
\end{align*}

It remains to integrate $\chi_M(lm)$ as a character in $M$, so that we can recognize a genuine double Dirichlet series. In order to do so, we can sum over classes of $lm$ in $H_C$, which is done at a cost of a uniformly finite linear combination by Lemma \ref{lem:sieving-classes-Hc}. We are therefore reduced to sums of the type
\begin{align*}
  \sum_{l|r/f_r(\beta)} \frac{\mu \beta(l)}{|l|^{w}}  &\sum_{m | r/f_r(\beta)} \frac{\beta(m)}{|m|^{1-w} } \frac{\varepsilon(w, \chi_E\beta)}{|E_0|^{\frac{1}{2}-w} \beta(E_0)} \prod_{\nu \in S} \frac{L_\nu(1-w, \chi_E \beta)}{L_\nu(w, \chi_E \beta)}  \rho(lm) Z^{S_r}\left(s+w-\tfrac{1}{2}, 1-w; \pi, \alpha \beta  \chi_{lm}\rho', \beta\right),
\end{align*}

where $\rho, \rho' \in \widehat{H}_C$, giving the claimed functional equation in $w$. \qed

After taking the sum over the ray classes $E$, we get an explicit functional equation relating $Z^{S_r}(s, w; \pi, \alpha, \beta)$ and $Z^{S_r}(\psi(s,w); \pi, \beta)$. In total, we get the relation between the double Dirichlet series at $(s,w)$ under the two transformations
\begin{align*}
\phi(s, w) & = \left(1-s, w+3s-\tfrac{3}{2}\right), \\
\psi(s, w) & = \left(s+w-\tfrac{1}{2}, 1-w\right).
\end{align*}

In particular, we can verify that the two involutions $\phi$ and $\psi$ generate a group of functional equations isomorphic to the dihedral group $D_6$. Indeed, it admits the presentation $\phi^2 = \psi^2 =1$ and $(\phi\psi)^6=1$. 

\section{Meromorphic continuation}
\label{sec:meromorphic-continuation}

We now prove that, based on the functional equations satisfied by the corrected double Dirichlet series, established in the previous section, we can extend $Z^S(s, w; \pi, \alpha, \beta)$ meromorphically to the whole $\C^2$. The method is essentially the same as in \cite{bump_sums_2004}. We will repeatedly use the following continuation principle due to Hartog \cite{hormander_introduction_1990}. 

\begin{thm}[Hartog's continuation principle]
Let $R$ be a connected tube domain, that is to say a connected domain of the form $S(\omega) = \{s \in \C^2 \ : \ \Re(s) \in \omega\}$ where $\omega$ is an open set of $\R^2$. Then any holomorphic function on $S(\omega)$ can be analytically continued to its convex hull $S(\widehat{\omega})$.
\end{thm}

\begin{prop}[Rough meromorphic continuation]
\label{prop:meromorphic-continuation}
Let $\alpha, \beta$ be characters of finite order, and $E, E'$ two classes in $H_C$. The function
\begin{equation}
(w-1) Z^S(s, w; \pi, \alpha \delta_E, \beta \delta_{E'}) 
\end{equation}

\noindent has an analytic continuation to the region $R_{1}$ made of the $(s,w) \in \C^2$ in
\begin{equation*}
\left\{
\begin{array}{cl}
\displaystyle \tau > \tfrac{5}{2} - 3\sigma & \text{if} \quad \sigma \leqslant -\frac{5}{14}- \varepsilon \\[1em]
\displaystyle  \tau > \tfrac{85}{28} - \tfrac{3}{2}\sigma & \text{if} \quad -\frac{5}{14}-\varepsilon < \sigma < \frac{19}{14} + \varepsilon \\[1em]
\tau > 1 & \text{if} \quad \sigma \geqslant \frac{19}{14}+\varepsilon
\end{array}
\right\}
\cup
\left\{
\begin{array}{cl}
\displaystyle  \sigma > \tfrac{13}{7} - \tau & \text{if} \quad \tau \leqslant - \varepsilon \\[1em]
\displaystyle  \sigma > \tfrac{13}{7} - \tfrac{1}{2}\tau & \text{if} \quad -\varepsilon < \tau < 1 + \varepsilon \\[1em]
\sigma > \tfrac{19}{14} & \text{if} \quad \tau \geqslant 1+\varepsilon
\end{array}
\right\}.
\end{equation*}
\end{prop}

\proof  First of all, let us deal with the functional equation corresponding to $\phi$. By definition, we have
\begin{equation}
Z^S(s, w; \pi, \alpha, \beta \delta_{E'}) = \sum_{D \in [E']} \frac{L^S(s, \pi \otimes \chi_D \alpha)}{|D|^w} \beta(D) a^S(s, D, \pi, \alpha).
\end{equation}

On $\sigma \geq \frac{19}{14}+\varepsilon$, the $L$-functions $L(s, \pi \otimes \chi_D \alpha)$ are uniformly convergent and therefore uniformly bounded. Moreover, $\beta$ happens to be a finite order character, so that it has norm less than one. Finally, the correction factor satisfies $a^S(s, D, \pi, \alpha) \ll |D|^\varepsilon$ for $\sigma > 1+\varepsilon$ by \eqref{correction-factor-bound-a}. So that the double Dirichlet series $Z^S(s, w; \pi, \alpha, \beta \delta_{E'})$ converges on $\tau > 1 + \varepsilon$, and is therefore holomorphic on the region $\sigma > \frac{19}{14}, \ \tau >1$. 

Appealing to the functional equation \eqref{fe1}, since all the extra factors appearing are holomorphic, we deduce that $Z^S(s, w; \pi, \alpha, \beta)$ is also holomorphic on $\sigma < -\frac{5}{14}$ and $\tau + 3\sigma > 5/2$. In the remaining region $-\frac{5}{14} \leqslant \sigma \leqslant \frac{19}{14}$, the bound $L(s, \pi \otimes \chi_D \alpha) a^S(s,D,\pi,\alpha) \ll |D|^{3/2 - 3\sigma}$ for $\sigma < -\frac{5}{14}$ and the Phragmén-Lindelöf principle yield that $L^S(s, \pi \otimes\chi_D \alpha) a^S(s,D,\pi,\alpha)$ is bounded by $|D|^{57/28 - 3\sigma/2}$. In particular, for $\tau > 85/28 - 3\sigma/2+\varepsilon$, the series $Z^S(s, w; \pi, \alpha, \beta \delta_E)$ uniformly converges and is therefore holomorphic in this region.

This proves that $Z^S(s, w; \pi, \alpha, \beta \delta_{E'})$ is holomorphic in the region $R_{1,1}$ made of all the $(s, w)$ in $\C^2$ such that
\begin{equation}
\left\{
\begin{array}{cl}
\displaystyle \tau > \tfrac{5}{2} - 3\sigma & \text{if} \quad \sigma \leqslant - \frac{5}{14} - \varepsilon, \\[1em]
\displaystyle  \tau > \tfrac{85}{28} - \tfrac{3}{2}\sigma & \text{if} \quad -\frac{5}{14} - \varepsilon < \sigma < \frac{19}{14} + \varepsilon, \\[1em]
\tau > 1 & \text{if} \quad \sigma \geqslant \frac{19}{14}+\varepsilon.
\end{array}
\right.
\end{equation}

Now, we turn to the formulation of the double Dirichlet series $Z^S(s, w; \pi, \alpha, \beta \delta_{E'})$ in terms of $\GL(1)$ $L$-functions given by \eqref{basic-identity}, that is
\begin{equation}
Z^S(s, w; \pi, \alpha \delta_E, \beta) = \sum_{M \in [E]} \frac{L^S(w, \chi_M \beta) }{|M|^s }b^S(w, M, \pi, \beta) \alpha(M).
\end{equation}

The completing factor satisfies $b^S(w, M, \pi, \beta) \ll a_\pi(M)$ for $\tau > 1+\varepsilon$ by \eqref{correction-factor-bound-b} and $\alpha$ is of norm bounded by one since it is a finite order character. In the domain $\tau > 1+\varepsilon$, the $L$-function $L^S(w, \chi_M \beta)$ is uniformly convergent and therefore bounded so that the double Dirichlet series above converges for $\sigma > \frac{19}{14}$. In particular, it is holomorphic on the region $\sigma >\frac{19}{14}, \ \tau > 1$. 

If $\tau \leqslant - \varepsilon$, we appeal to the functional equation \eqref{fe-w} and note that all the extra factors appearing there are holomorphic. This shows $Z^S(s, w; \pi, \alpha, \beta)$ is holomorphic on $\tau < 0$, $\sigma > \frac{13}{7} - \tau$. In order to interpolate these bounds in between by the Phragmén-Lindelöf principle, it is necessary to get rid of the simple pole of $L^S(w, \chi_M \beta)$ at $w=1$ when the character is trivial, which is achieved by multiplying the factor $w-1$. The Phragmén-Lindelöf principle therefore gives $(w-1) L^S(w, \chi_M \beta) b^S(w,M,\pi,\beta) \ll |M|^{6/7 -\tau/2}$, so that the series converges in the region $\sigma > 13/7 - \tau/2 +\varepsilon$. Altogether, we conclude that $(w-1)Z^S(s, w; \pi, \alpha \delta_E, \beta)$ is holomorphic in the region $R_{1,2} \subseteq \C^2$ made of $(s, w) \in \C^2$ such that
\begin{equation}
\label{region-R1-1}
\left\{
\begin{array}{cl}
\displaystyle \sigma > \tfrac{13}{7} - \tau & \text{if} \quad \tau \leqslant - \varepsilon, \\[1em]
\displaystyle \sigma > \tfrac{13}{7} - \tfrac{1}{2}\tau & \text{if} \quad -\varepsilon < \tau < 1 + \varepsilon, \\[1em]
\sigma > \tfrac{19}{14} & \text{if} \quad \tau \geqslant 1+\varepsilon,
\end{array}
\right.
\end{equation}

\noindent ending the proof with $R_1 = R_{1,1} \cup R_{1,2}$. \qed

We apply now repeatedly the functional equations and the Hartog continuation principle to get the meromorphic continuation on the whole $\C^2$. Introduce the functions
\begin{align*}
\Phi(s, w) & = \prod_{P \in S_f} \prod_j \left(1- \frac{\alpha \gamma_j(P)^2}{ |P|^{2-2s}} \right), \\
\Psi(s,w) & =  \prod_{P \in S_f}  \left(1-\frac{\beta(P)^2}{|P|^{2-2w}} \right), \\
P(s, w) & = w(w-1)\left(w+3s-\tfrac{3}{2}\right)\left(w+3s-\tfrac{5}{2}\right)(3s+2w-3).
\end{align*}

\begin{prop}[Meromorphic continuation]
\label{prop:meromorphic-continuation-old}
Let 
\begin{equation*}
\xi(s, w) = P(s,w)\Phi(s, w) \Phi(\phi(s, w)) \Phi(\psi(s, w)) \Phi(\psi \phi(s, w)) \Psi(s, w) \Psi(\phi(s, w)).
\end{equation*}
Then, the completed double Dirichlet series
\begin{equation}
\xi(s, w) Z^S(s, w; \pi, \alpha, \beta)
\end{equation}
admits an analytic continuation to $\C^2$.
\end{prop}

\proof Up to a uniformly finite linear combination, Lemma \ref{lem:sieving-classes-Hc} ensures that we are reduced to prove that, for every classes $E, E'$ in $H_C$, the completed function $\xi(s, w) Z^S(s, w; \pi, \alpha \delta_E, \beta \delta_{E'})$ admits an analytic continuation to $\C^2$. We use the functional equations combined with the Hartog continuation principle to extend the previous domain of holomorphy $R_1$ of $(w-1)Z^S(s,w; \pi, \alpha \delta_E, \beta \delta_{E'})$ to the whole $\C^2$. Indeed, applying the functional equation \eqref{fe1} transforming by $\phi$, we get that 
\begin{equation}
\Phi(s, w) \left(w+3s-\tfrac{5}{2}\right) Z^S(s,w; \pi, \alpha \delta_E, \beta\delta_{E'})
\end{equation}
is analytic on $\phi(R_1)$, so that adding all the completing factors we deduce that 
\begin{equation}
\Phi(s, w) (w-1)\left(w+3s-\tfrac{5}{2}\right) Z^S(s,w; \pi, \alpha \delta_E, \beta\delta_{E'})
\end{equation}
is analytic on $R_1 \cup \phi(R_1)$. The Hartog continuation principle therefore allows to analytically continue it to the convex hull $R_2$ of the union $R_1 \cup \phi(R_1)$. 

Applying now the functional equation \eqref{fe-w} corresponding to $\psi$, we get that
\begin{align*}
& \Psi(s, w) \Phi(\psi(s, w)) w (3s+2w-3) Z^S(s,w; \pi, \alpha \delta_E, \beta\delta_{E'})
\end{align*}

\noindent is analytic on $\psi(R_2)$, so that adding all the completing factors we deduce that 
\begin{align*}
& \Psi(s, w) \Phi(s, w) \Phi(\psi(s, w)) w (w-1)\left(w+3s-\tfrac{5}{2}\right) (3s+2w-3) Z^S(s,w; \pi, \alpha \delta_E, \beta\delta_{E'})
\end{align*}

is analytic on $R_2 \cup \psi(R_2)$. The Hartog continuation principle therefore allows to analytically continue it to the convex hull $R_3$ of the union $R_2 \cup \psi(R_2)$. 

Applying again the functional equation \eqref{fe1} transforming by $\phi$ and gathering all the extra factors, we get that $\xi(s, w) Z^S(s, w; \pi, \alpha\delta_E, \beta\delta_{E'})$ is analytic on the convex hull of $R_3 \cup \phi(R_3)$, that happens to be the whole $\C^2$. \qed 

\section{Sieving process}
\label{sec:sieving-process}

\subsection{Sieving out squares}

Now it remains to show that these good analytic properties of $Z^S(s, w; \pi, \alpha, \beta)$, in particular its analytic continuation proven in Proposition \ref{prop:meromorphic-continuation-old}, transfer into good analytic properties of the original double Dirichlet series without correction factors. 

Introduce the squarefree part of the pure double Dirichlet series, 
\begin{equation}
\label{Zstar-def}
Z^{S}_{\star}(s, w; \pi, \alpha, \beta) := \sum_{\substack{D \in I^+(S) \\ \text{squarefree}}} \frac{L^S(s, \pi \otimes \chi_D \alpha)}{|D|^w} \beta(D).
\end{equation}

\textit{Remark.} The correction factors play no role here since they are trivial on squarefree ideals by definition. The very reason of their introduction is precisely the existence of square parts that need to be compensated in order to write the functional equations.

Introduce, for $r \in I^+(S)$,
\begin{equation}
\label{Zr-def}
Z^S_r(s, w; \pi, \alpha, \beta)  :=  \sum_{\substack{D \in I^+(S) \\ D = D_0D_1^2 \\ r | D_1}} \frac{L^S(s, \pi \otimes \chi_D \alpha)}{|D|^w} a^S(s, D, \pi, \alpha) \beta(D).
\end{equation}

\begin{prop}
\label{prop:relation-Zstar-Zr}
For $\sigma, \tau > 1$, 
\begin{equation}
\label{relation-Zstar-Zr}
Z^{S}_{\star}(s, w; \pi, \alpha, \beta)  = \sum_{\substack{r \in I^+(S) }} \mu(r) Z_r^S(s, w; \pi, \alpha, \beta).
\end{equation}
\end{prop}

\proof By definition, the right hand side in \eqref{relation-Zstar-Zr} can be rewritten as
\begin{equation}
\label{ccc}
\sum_{\substack{r \in I(S) }} \mu(r) \sum_{\substack{D \in I(S) \\ D = D_0D_1^2 \\ r | D_1}} \frac{L^S(s, \pi \otimes \chi_D \alpha)}{|D|^w} a^S(s, D, \pi, \alpha) \beta(D).
\end{equation}

The proposition reduces to a property that holds for any arithmetic function $f(D)$ in general. Indeed, summing over ideals in $I^+(S)$, 
\begin{equation}
\sum_{r} \mu(r) \sum_{\substack{D \\ r | D_1}} f(D)  =  \sum_D f(D) \sum_{r | D_1} \mu(r) = \sum_{\substack{D \\ \text{squarefree}}} f(D) ,
\end{equation}

\noindent by Möbius inversion, since the sum of $\mu(r)$ is equivalent to selecting those $D$ for which $D_1=1$, that is to say squarefree ideals. Taking $f(D)$ to be the summands in \eqref{ccc} yields the result. \qed

According to \eqref{relation-Zstar-Zr}, we are reduced to studying double Dirichlet series of type $Z_r^S(s, w; \pi, \alpha, \beta)$ for squarefree ideals $r$. Let $l \in I^+(S)$. Introduce 
\begin{equation}
\label{Z$L$-def}
Z^{S}_{(l)}(s, w; \pi, \alpha, \beta) := \sum_{\substack{D \in I^+(S) \\ D = D_0D_1^2 \\ (D_1, l) = 1}} \frac{L^S(s, \pi \otimes \chi_D \alpha)}{|D|^w} a^S(s, D, \pi, \alpha) \beta(D).
\end{equation}
\begin{lem}
\label{lem:expression-Zr}
For $\sigma, \tau > 1$, and $r \in I^+(S)$ squarefree,
\begin{equation}
Z^S_r(s, w; \pi, \alpha, \beta)  =  \sum_{l | r} \mu(l) Z^{S}_{(l)}(s, w; \pi, \alpha, \beta).
\end{equation}
\end{lem}

\proof This relation holds more generally for any arithmetic function $f(D)$. In particular, by switching summations, we get
\begin{equation}
\sum_{l | r} \mu(l) \sum_{\substack{D \\ (D_1, l) = 1}} f(D) = \sum_{D} f(D) \sum_{l | \frac{r}{(r,D_1)}} \mu(l).
\end{equation}

This last relation is justified since the relations $l | r$ and $(D_1, l)=1$ rewrite $r = (r, D_1) l k$ for an integer ideal $k$, that is to say $l$ divides $r/(r,D_1)$. The result follows by Möbius inversion, for the last sum translates into the condition $r=(r, D_1)$, that is to say $r | D_1$. With $f(D)$ the summand in \eqref{Z$L$-def}, this is the claimed result. \qed

This proposition along with the relation \eqref{relation-Zstar-Zr} shows that the squarefree pure double Dirichlet series $Z^{S}_{\star}(s, w; \pi, \alpha, \beta)$ can be expressed as a sum of the corrected double Dirichlet series $Z^S_{(l)}(s, w; \pi, \alpha, \beta)$. The following proposition relates the latter to the corrected double Dirichlet series $Z^{S_l}(s,w;\pi,\alpha,\beta)$. 

\begin{prop}
\label{prop:relation-Z$L$-Z(l)}
For every pair of characters $\alpha, \beta$ of finite order,  $l$ square-free,
\begin{align}
&  \prod_{P | l} \prod_j \left(1-\frac{\alpha\gamma_j(P)^2}{|P|^{2s}}\right)  Z^{S}_{(l)}(s, w; \pi, \alpha, \beta)
\end{align}

is a uniformly finite linear combination of expressions of the form
\begin{align}
\label{pp}
& \sum_{l_3 | l} \sum_{n_j | l_3} \frac{\mu \alpha^2\gamma_j^2(n_j)}{|l_3|^{w}|n_j|^{2s}} \sum_{m_j | l /l_3} \frac{\chi_{m_1m_2m_3} \beta \rho(l_3) \alpha \rho'(m_1m_2m_3) \gamma_j(m_j)}{|m_1m_2m_3|^s}  Z^{S_l} (s, w; \pi, \alpha \chi_{l_3}, \beta \rho \chi_{m_1m_2m_3}).
\end{align}
\end{prop}

\textit{Remark.} From \cite[Proposition 1.5]{li_determination_2007} we know that there is a suitable choice of representatives in Section \ref{subsec:quadratic-symbols} so that, for $D \in S$, the quadratic characters $\chi_D$ defined with respect to $S$ or to $S_l$ match.

\proof This is essentially \cite[Proposition 3.5]{chinta_determination_2005} or the analogous \cite[Proposition 4.14]{diaconu_multiple_2003} adapted to the number field case. Since none of these articles prove the result explicitly, we provide the details here. 

By definition, we have
\begin{equation}
\label{1}
Z^{S}_{(l)}(s, w; \pi, \alpha, \beta) = \sum_{(D_1, l)=1} \frac{L^S(s, \pi \otimes \chi_D \alpha)}{|D|^w}\beta(D)a^S(s, D, \pi, \alpha).
\end{equation}

Introduce the variable $l_3 = (D_0, l)$, so that $l_3$ is squarefree as $D_0$. Making the change of variables $D_0 \to D_0l_3$ and $l \to ll_3$ yields the restriction $(D_0,l) = 1$, which combined with the condition on $D_1$ gives $(D,l) = 1$. We can rewrite the above as
\begin{equation}
Z^{S}_{(l)}(s, w; \pi, \alpha, \beta) = \sum_{\substack{l_3 | l}} \sum_{(D, l)=1} \frac{L^S(s, \pi \otimes \chi_{Dl_3} \alpha)}{|D|^w|l_3|^w}\beta(Dl_3)a^S(s, Dl_3, \pi, \alpha).
\end{equation}

Making appear explicitly the local $L$-factor relative to the place $l$ yields
\begin{equation}
L^S(s, \pi \otimes \chi_{Dl_3}\alpha) = L^{S_l}(s, \pi \otimes \chi_{Dl_3}\alpha) \prod_{P | l}  \prod_{j=1}^3  \left( 1-\frac{\chi_{Dl_3} \gamma_j(P)}{|P|^{s}} \right)^{-1}.
\end{equation}

Introducing this expression in \eqref{1} and multiplying both sides by $\prod_{P | l} \prod_j (1-\alpha\gamma_j(P)^2|P|^{-2s})$ yields
\begin{align*}
& \prod_{i=1}^3 \prod_{P | l} \left(1-\frac{\alpha\gamma_j(P)^2}{|P|^{2s}}\right) Z^{S}_{(l)}(s, w; \pi, \alpha, \beta) \\
& \qquad  = \sum_{\substack{l_3}} |l_3|^{-w} \prod_{P | l} \prod_j \left(1-\alpha \gamma_j(P)^2|P|^{-2s}\right) \sum_{(D, l)=1} \frac{L^{S_l}(s, \pi \otimes \chi_{Dl_3}\alpha) }{|D|^w}\\
& \qquad \qquad \times \prod_{P | l} \prod_j \left( 1-\frac{\alpha\chi_{Dl_3} \gamma_j(P)}{|P|^{s}} \right)^{-1}  \beta(Dl_3)a^{S_l}(s, D, \pi, \alpha \chi_{l_3}).
\end{align*}
We have utilized \eqref{extraprop} to get the final form of the correcting factor.

By definition of the quadratic symbol, $\chi_{Dl_3}(P)$ is zero as soon as $P | Dl_3$, that is to say $P | l_3$ since $(D, l)=1$. Note also that, since $l_3$ is squarefree, the conditions $P | l$ and  $P \nmid l_3$ can be summarized as $P|l/l_3$. Moreover, since the characters $\chi_D$ are quadratic, we have for all $j$
\begin{align*}
\prod_{P \nmid l_3} \frac{1 - \alpha\gamma_j(P)^2|P|^{-2s}}{1-\alpha \chi_{Dl_3} \gamma_j(P)|P|^{-s}} &= \prod_{P \nmid l_3} \left(1+\frac{\alpha \chi_{Dl_3}\gamma_j (P)}{ |P|^{s}}\right).
\end{align*}

By developing the involved Euler product, we get
\begin{equation}
\prod_{P \nmid l_3} \left(1+\frac{\alpha\chi_{Dl_3} \gamma_j (P)}{ |P|^{s}}\right) = \sum_{m | l/l_3} \frac{\alpha \chi_{Dl_3}\gamma_j(m)}{ |m|^{s}}.
\end{equation}

Altogether, the expression above rewrites
\begin{align*}
& \prod_{j=1}^3 \prod_{P | l} \left(1-\frac{\alpha \gamma_j(P)^2}{ |P|^{2s}}\right) Z^{S}_{(l)}(s, w; \pi, \alpha, \beta) \\
=& \sum_{\substack{l_3 | l}} |l_3|^{-w} \prod_{P | l_3} \prod_j  \left(1-\frac{\alpha \gamma_j(P)^2}{ |P|^{2s}}\right) \sum_{(D, l)=1} \frac{L^{S_l}(s, \pi \otimes \chi_{Dl_3}\alpha) }{|D|^w}  \sum_{m_j | l/l_3} \frac{\alpha \chi_{Dl_3} (m_1m_2m_3) \gamma_j(m_j)}{|m_1m_2m_3|^s}  \beta(Dl_3)a^{S_l}(s, D, \pi, \alpha \chi_{l_3}).
\end{align*}

By the quadratic reciprocity law \eqref{quadratic-reciprocity-law}, we get
\begin{equation}
\chi_{Dl_3}(m_1m_2m_3) = \eta(Dl_3, m_1m_2m_3) \chi_{m_1m_2m_3}(Dl_3) .
\end{equation}

By the orthogonality relations \eqref{orthogonality-relation} we can select the ray classes in $H_C$ for $Dl_3$ as well as $m_1m_2m_3$ at the cost of a uniformly finite linear combination, and be reduced to sums of the form
\begin{align*}
& \sum_{l_3 | l} \sum_{n_j | l_3} \frac{\mu \alpha^2\gamma_j^2(n_j)}{|l_3|^{w}|n_j|^{2s}} \sum_{m_j | l /l_3} \frac{\chi_{m_1m_2m_3} \beta \rho(l_3) \alpha \rho'(m_1m_2m_3) \gamma_j(m_j)}{|m_1m_2m_3|^s}  Z^{S_l} (s, w; \pi, \alpha \chi_{l_3}, \beta \rho \chi_{m_1m_2m_3}),
\end{align*}
where $\rho, \rho' \in \widehat{H}_C$. This is the claimed result. \qed

\subsection{Vertical bounds}
\label{subsec:vertica$L$-bounds}

The relation \eqref{relation-Zstar-Zr} is unfortunately not a linear combination, but an infinite sum. As such, extra bounds in the $r$-aspect are necessary in order to ensure convergence so that the sum remains meromorphic. This is the aim of this subsection and also precisely where Chinta and Diaconu were not able to generalize their results to all number fields, due to the lack of large sieve inequalities in this setting. They relied on a quadratic large sieve result due to Heath-Brown \cite{heath-brown_mean_1995}. Recent results due to Goldmakher and Louvel \cite{goldmakher_quadratic_2013} open the path beyond this point.

\begin{lem}[Goldmakher-Louvel]
\label{lem:GL}
Let $F$ be a number field and $X$ be the set of quadratic characters $\chi_\d$ with $\d$ integral ideals out of $S$. Let $(\lambda_\d)_\d$ be a sequence of complex numbers parametrized by integral ideals $\d$ of $F$. Then we have for every $\varepsilon > 0$ and every $M, N \geqslant 1$,
\begin{equation}
\sum_{|A| \leqslant M}^\star \left| \sum_{|D| \leqslant N}^\star \lambda_\d \chi_\d(A) \right|^2 \ll_{\varepsilon} (MN)^\varepsilon (M+N) \sum_{|\d|  \leqslant N}^\star |\lambda_\d|^2	,
\end{equation}

\noindent where the starred sums stands for sums restricted to squarefree ideals.
\end{lem}

We deduce from this quadratic large sieve inequality for general number fields the following estimate on mean square of twisted central values.

\begin{lem}
\label{lem:convergence-$L$-sum}
Let $\pi$ be a self-contragredient cuspidal automorphic representation of $\GL(3)$ over $F$. For all $\varepsilon > 0$ and character $\alpha$ of finite order we have, for all $Y>0$,
\begin{equation}
\sum_{|D| \leqslant Y} \left| L\left( \tfrac{1}{2}, \pi \otimes \chi_\d \alpha \right) \right|^2 \ll_\varepsilon Y^{3/2 + \varepsilon}.
\end{equation}
\end{lem}

%
%

For the remainder of the section the bounds for $Z^{S_r}\left(s, w ; \pi,\alpha,  \beta\right)$ are given for $s=\frac{1}{2}$. In practice, we will need them around a small neighborhood of $s=\frac{1}{2}$, and it is straightforward to check that the bounds below still hold for $s$ in a small compact neighborhood of $\frac{1}{2}$. 

\begin{prop}
\label{prop:bounds-r-aspect-right}
For $\tau > 5/4 + \varepsilon$ and characters $\alpha, \beta$ of finite order, 
\begin{equation}
\label{bound-right}
Z^{S_r}\left(\tfrac{1}{2}, w ; \pi,\alpha,  \beta\right) \ll_\varepsilon |r|^\varepsilon.
\end{equation}
\end{prop}

\proof Proposition \ref{prop:meromorphic-continuation-old} states that the function $Z^{S_r}\left(\frac{1}{2}, w ; \pi,\alpha,  \beta\right)$ is analytic in $w$ except for possible poles at $w=0$, $w=\frac{3}{4}$ and $w=1$. $Z^{S_r}\left(\frac{1}{2}, w; \pi, \alpha, \beta\right)$ converges as a consequence of  Lemma \ref{lem:convergence-$L$-sum}. We can be see this by writing down partial sums and recalling $\chi_D$ only depends on the squarefree part $D_0$ of $D$, getting
\begin{align*}
& \sum_{\substack{|D| \leqslant Y \\ (\d,r) = 1}} \frac{L^{S_r}\left(\frac{1}{2}, \pi , \chi_\d \alpha\right)}{|D|^w} \beta(D) a^{S_r}\left(\frac{1}{2}, D, \pi, \alpha\right)\\
& \qquad  \ll_\varepsilon |r|^\varepsilon \sum_{\substack{|D_0| \leqslant Y \\ (\d_0, r)=1}} \frac{\left| L^{S_r}\left( \frac{1}{2}, \pi \otimes \chi_{\d_0}\alpha \right) \right|}{|D_0|^\tau} \sum_{\substack{|D_1|^2 \leqslant Y/|D_0| \\ (\d_1, r) = 1}}  \frac{\left| a^{S_r}\left(\frac{1}{2}, D, \pi, \alpha\right) \right|}{|D_1|^{2\tau}} \\
& \qquad \ll_\varepsilon |r|^\varepsilon \sum_{\substack{|D_0| \leqslant Y \\ (\d_0, r)=1}} \frac{\left| L^{S_r}\left( \frac{1}{2}, \pi \otimes \chi_{\d_0} \alpha \right) \right|}{|D_0|^\tau} \sum_{ (\d_1, r) = 1}  \frac{\left| a^{S_r}\left(\frac{1}{2}, D, \pi, \alpha\right)  \right|}{|D_1|^{2\tau}} \\
& \qquad \ll_\varepsilon |r|^\varepsilon \sum_{(\d_0, r)=1} \frac{\left| L^{S_r}\left( \frac{1}{2}, \pi \otimes \chi_{\d_0}  \alpha  \right) \right|}{|D_0|^\tau} 
\end{align*}

\noindent and this last sum is convergent  for $\tau > 5/4$ by Lemma \ref{lem:convergence-$L$-sum}. This gives the claim. \qed

\begin{prop}
\label{prop:bounds-r-aspect-left}
For $\tau = -1/4-\varepsilon$ and $\alpha, \beta \in \widehat{H}_C$,
\begin{equation}
Z^{S_l}\left(\frac{1}{2}, w ; \pi, \alpha, \beta\right) \ll_{\varepsilon} |l|^{5+\varepsilon} \sum_{\rho \in \widehat{H}_C} \sum_{D_0} \frac{\left| L\left( \frac{1}{2}, \pi \otimes \chi_D \rho \right) \right|}{|D_0|^{5/4+\varepsilon} |C_\rho|^{1/4}}.
\end{equation}
\end{prop}

\proof This is analogous to \cite[Proposition 3.4]{chinta_determination_2005}. First of all, we know explicit bounds on the local $L$-factors as a consequence of the Stirling formula, see for instance \cite[Equation (3.5)]{bump_sums_2004}. We have for all $\sigma_1 > \sigma_2$ and large enough $|t|$, and for $\pi$ an automorphic cuspidal representation on $\GL_n$ for $n \in \{1, 3\}$.
\begin{equation}
\prod_{v \in S_\infty} \frac{L_v(\sigma_1+ it, \pi \otimes \chi_D)}{L_v(\sigma_2 - it, \pi \otimes \chi_D)} \ll (|t|+1)^{n (\sigma_1 - \sigma_2)/2}.
\end{equation}

Moreover, since every character appearing in the functional equations is of finite order and therefore of norm bounded by one, and the functional equations involve uniformly finite linear combinations, \eqref{fe1} and \eqref{fe-w} imply the following bounds.
\begin{lem}
\label{lem:fe-bounds}
For an automorphic cuspidal representation $\pi$ on $\GL(3)$, $\alpha, \beta \in \widehat{H}_C$, $E$ a class in $H_C$ prime to $\, S$, we have $Z^{S_r}(s, w; \pi, \alpha, \beta \delta_E)$ bounded by a uniformly finite linear combination of expressions, for $\rho \in \widehat{H}_C$,
\begin{align}
\label{fe-s-bound}Z^{S_r}(s, w; \pi, \alpha, \beta \delta_E) & \ll f(\alpha)^{3(\frac{1}{2}-\sigma)}\prod_{P | r / f_r(\alpha)} \prod_j \left|1- \frac{ \alpha \gamma_j(P)^2}{|P|^{2-2s}}  \right| ^{-1} \sum_{l_j | r / f_r(\alpha)} \frac{ |\gamma_j(l_j)|}{|l_j|^{1-\sigma}} \\
\notag& \qquad \times \sum_{m_j | r/f_r(\alpha)}  \frac{|\gamma_j(m_j)|}{|m_j|^{\sigma}}     \left|Z^{S_r}(\phi(s, w) ; \pi, \alpha , \beta \rho \chi_\pi \chi_{l_1l_2l_3m_1m_2m_3})\right|,
\end{align}
and $Z^{S_r}(s, w; \pi, \alpha \delta_E, \beta)$ bounded by a uniformly finite linear combination of expressions, for $\rho \in \widehat{H}_C$,
\begin{align}
\label{fe-w-bound} Z^{S_r}(s, w; \pi, \alpha \delta_E, \beta) & \ll  f(\beta)^{\frac{1}{2}-\nu} \!\! \prod_{P|r/f_r(\beta)}  \left|1-\frac{\beta(P)^2}{|P|^{2-2w}}\right|^{-1} \!\!\! \sum_{l|r/f_r(\beta)} \frac{1}{|l|^{w}} \sum_{m | r/f_r(\beta)} \frac{1}{|m|^{1-w}}   \left|Z^{S_r}(\psi(s, w); \pi, \alpha \rho \beta \chi_{lm}, \beta)\right|.
\end{align}
\end{lem}

Note that $\psi\phi\psi\phi\psi$ maps the line $(\frac{1}{2}, \frac{5}{4}+\varepsilon+it)$ onto the line $(\frac{1}{2}, -\frac{1}{4}-\varepsilon-it)$. Because of the bound \eqref{bound-right} on the former, the strategy is to apply successively the functional equations \eqref{fe1} and \eqref{fe-w} in order to apply the bounds on the line $(\frac{1}{2}, \frac{5}{4}+\varepsilon+it)$, carefully gathering and bounding all the extra terms appearing in the bounds provided by Lemma \ref{lem:fe-bounds}.

The elegant matrix formulation of \cite[pages 2951-2955]{chinta_determination_2005} stands \textit{mutatis mutandis} in our case. Indeed, the bound obtained in Proposition \eqref{fe-s-bound} is analogous to their bound (3.16), and the bound obtained in Proposition \eqref{fe-w-bound} is analogous to their bound (3.17). Since the exponents are the same everywhere, we do not repeat their computations except for the following fact.

 A crucial input here is the bound towards Ramanujan $|\alpha_\pi(P)|^{-1} \ll |P|^{5/14}$ due to \cite{blomer_ramanujan_2011} for general number fields. Indeed, using only $|\alpha_\pi(P)|^{-1} \ll |P|^{1/2}$, the matrix bound \cite[(3.31)]{chinta_determination_2005} and computing the explicit dependency in $a_\pi(P)$ gives an entry-by-entry bound by
\begin{equation}
\left( 
\begin{matrix}
|P|^5 & |P|^{19/4} & |P|^5 & 0 \\
|P|^{19/4} & |P|^{9/2} & |P|^{19/4} & 0 \\
|P|^5 & |P|^{19/4} & |a_\pi(P)|^{2} |P|^4 & 0 \\
0& 0& 0& |P|^{15/2}
\end{matrix}
\right).
\end{equation}

Inspecting in particular the $(3,3)$-entry, it is less than $|P|^{19/4}$ (the bound used by Chinta and Diaconu) if and only if $|a_\pi(P)| \ll |P|^{3/8}$. In particular, it is true by the Blomer-Brumley bound. Hence the final bound on $Z^{S_r}(s, w; \pi, \alpha \delta_E, \beta)$ given in \cite[(3.32)]{chinta_determination_2005} remains valid for every cuspidal automorphic representation of $\mathrm{GL}(3)$. The remaining of the proof is straightforward. \qed

\subsection{Analytic continuation}

We have now the tools required in order to analytically continue the Dirichlet series $Z^S_\star(s, w; \pi, \alpha, \beta)$. 

\begin{prop}
\label{prop:meromorphic-continuation-Zstar}
The series $Z^S_\star(s, w; \pi, \alpha, \beta)$ has meromorphic continuation to a tube domain containing $(\frac{1}{2}, 1)$. More precisely, the function $\xi(s,w)Z^S_\star(s, w; \pi, \alpha, \beta)$ has analytic continuation in the region of the $(s, w) \in \C^2$ defined by 
\begin{equation}
\label{region}
\left\{ \sigma > \tfrac{1}{2}, \tau >  \tfrac{119}{124} \right\} \bigcup \left\{0 \leqslant \sigma \leqslant \tfrac{1}{2}, \tau >  -\tfrac{191}{62} \sigma + \tfrac{5}{2} \right\}.
\end{equation}
\end{prop}

\proof  We follow closely \cite[Proposition 3.1]{chinta_determination_2005} here. Indeed, we obtained statements analogous to the ones used in their proof so that we are now able to conclude in the same way. We include the relevant details in our setting for completeness.

By Proposition \ref{prop:relation-Zstar-Zr}, $Z^S_\star(s, w; \pi, \alpha, \beta)$ is a sum over $r\in I^+(S)$ of functions $Z^S_r(s, w; \pi, \alpha, \beta)$, so it is enough to establish bounds on this last quantity in the $r$-aspect. We need an estimate in the central strip between $-1/4$ and $5/4$.

Let $\sigma \geqslant \tfrac{1}{2}$ and $\tau > \tfrac{5}{4}$. Recall that, by definition, 
\begin{equation}
\label{Zr-recall}
Z^S_r(s, w; \pi, \alpha, \beta)  =  \sum_{\substack{D \in I^+(S) \\ D = D_0D_1^2 \\ r | D_1}} \frac{L^S(s, \pi \otimes \chi_D \alpha)}{|D|^w} a^S(s, D, \pi, \alpha) \beta(D).
\end{equation}

The explicit description of the correction factors in \cite{bump_sums_2004} yields the bounds, see for instance the paper \cite[(3.39)]{chinta_determination_2005} (correcting an obvious typo in the denominator),
\begin{equation}
\sum_{D_1 \in I^+(S)} \frac{a^S(s, r^2 D, \pi, \alpha)}{|rD_1|^{2w}}  \ll |r|^{-2\tau} \left( (|b_\pi(r)| + |a_\pi(r)|)|r|^{1-2\sigma} + 1 \right), 
\end{equation}

\noindent where $b_\pi(r)$ are the coefficients of $L(s, \pi, \mathrm{sym}^2)$ and the bound is uniform in $D_0$. Plugging this bound into \eqref{Zr-recall} and taking $\sigma \geqslant \tfrac{1}{2}$ and $\tau > \tfrac{5}{4}$, a region in which the sum over $D_0$ converges, we get 
\begin{equation}
\label{BB}
Z^S_r(s, w; \pi, \alpha, \beta) \ll |r|^{-2\tau} \left( (|b_\pi(r)| + |a_\pi(r)|)|r|^{1-2\sigma} + 1 \right).
\end{equation}

We now turn to bounds on the domain $\tau < -1/4$. By Lemma \ref{lem:expression-Zr}, $Z^S_r(s, w; \pi, \alpha, \beta)$ is a finite sum of $Z^S_{(l)}(s, w; \pi, \alpha, \beta)$, almost uniformly finite (the sum over $l \mid r$ will amount at most to a factor $|r|^\varepsilon$ which is not impacting the convergence domains). It is therefore enough to bound $Z^S_{(l)}(s, w; \pi, \alpha, \beta)$,  for $l \mid r$, in the $r$-aspect.

By Proposition \ref{prop:relation-Z$L$-Z(l)}, $Z^S_{(l)}(s, w; \pi, \alpha, \beta)$ is in turn a finite sum of corrected double Dirichlet series $Z^{S_l}(s, w; \pi, \alpha, \beta)$. For these, Proposition \ref{prop:bounds-r-aspect-left} provides explicit bounds on the vertical line $\tau = -\tfrac{1}{4} - \varepsilon$. Plugging these bounds in the expression \eqref{pp} and noting that all the sums in \eqref{pp} are uniformly bounded in $r$ except the sum of $|l_3|^{-w}$ over $l_3$, that contributes as $r^{1/4}$, we deduce that $Z^S_{(l)}(s, w; \pi, \alpha, \beta)$ is a uniformly finite linear combination of expressions bounded by
\begin{equation}
\label{AA}
|r|^{21/4 + \varepsilon} \sum_{\rho \in \widehat{H}_{Cl}} \sum_{D_0 \ \text{squarefree}} \frac{|L(\tfrac{1}{2}, \pi \otimes \chi_{D_0} \rho)|}{|D_0|^{5/4}|C_\rho|^{1/4}}, 
\end{equation}

\noindent where $C_\rho$ is the conductor of $\rho$. The completed $\xi(s,w) Z^S_{(l)}(s, w; \pi, \alpha, \beta)$ is holomorphic on $\C^2$ by Proposition \ref{prop:meromorphic-continuation-old}, so that we can apply the convexity principle to the two vertical bounds \eqref{BB} on $\tau = \tfrac{5}{4}+\varepsilon$ and \eqref{AA} on $\tau = -\tfrac{1}{4} - \varepsilon$.  We get, for $-\tfrac14 - \varepsilon < \tau < \tfrac54 + \varepsilon$, 
\begin{align}
\label{parenthesis}
Z^S_{(l)}(s, w; \pi, \alpha, \beta) & \ll |r|^{95/24 - 31\tau/6 + \varepsilon} \left( (|b_\pi(r)| + |a_\pi(r)|)|r|^{1-2\sigma} + 1 \right) \\
& \qquad \times \sum_{\rho \in \widehat{H}_C} \sum_{D_0 \ \text{squarefree}} \frac{|L(\tfrac{1}{2}, \pi \otimes \chi_{D_0} \rho)|}{|D_0|^{5/4}|C_\rho|^{1/4}}.
\end{align}

Putting together the summations appearing in \eqref{relation-Zstar-Zr} and \eqref{pp}, we hence deduce that $Z_\star(s, w; \pi, \alpha, \beta)$ is bounded by a sum of 3 terms, corresponding to the three terms appearing in the parenthesis \eqref{parenthesis}. The hardest to bound is the one containing the $b_\pi$ coefficients, given by
\begin{equation}
\sum_{r \ \text{squarefree}} |b_\pi(r)| \cdot  |r|^{95/24 - 31\tau/6 + \varepsilon} \sum_{\rho \in \widehat{H}_{Cl}} \sum_{D_0 \ \text{squarefree}} \frac{|L(\tfrac{1}{2}, \pi \otimes \chi_{D_0} \rho)|}{|D_0|^{5/4}|C_\rho|^{1/4}}.
\end{equation}

This sum converges as long as $\tfrac{31}{6}\tau - \tfrac{95}{24} >1$ exactly as in \cite[p. 2959]{chinta_determination_2005}. The crux of the argument is essentially that the Dirichlet series involving $\mathrm{sym}^2(\pi)$-coefficients converge in this case. Indeed, the sum analogous to \cite[(3.50) and (3.55)]{chinta_determination_2005} is convergent by \cite[Appendix 2]{kim_functoriality_2002}. This concludes the proof of the analytic continuation of  $Z_\star(s, w; \pi, \alpha, \beta)$ to 
\begin{equation}
\left\{ \sigma > \tfrac{1}{2}, \tau >  \tfrac{119}{124} \right\} .
\end{equation}

 Applying Hartog's principle between this domain and the domain $\sigma > 0$ and $\tau >~5/2$, obtained directly by the convexity bound $L(s, \pi \otimes \chi_D) \ll |D_0|^{3/2 + \varepsilon}$, yields the claimed result. \qed

\section{Residues and Fourier coefficients}
\label{sec:loca$L$-analysis}

\subsection{Analytic properties}
\label{subsec:loca$L$-analysis}

The previous section concluded with the meromorphic continuation of the double Dirichlet series $Z_\star(s, w; \pi, \alpha, \beta)$, around the point $(\frac{1}{2}, 1)$. Theorem \ref{thm:result} will eventually follow from a precise study of the corresponding residues at $(\frac{1}{2}, 1)$, which are related to the Fourier coefficients of $\pi$. It is necessary to write down more precisely the local behavior around the point $(\frac{1}{2}, 1)$.

Let $E$ be an element in $H_C$. Consider an ideal $r \in I^+(S)$ that is either $\O$ or a prime ideal. Introduce
\begin{equation}
Z_\star(s, w; \pi, 1, \delta_E\chi_r) := \sum_{\substack{D \in [E] \\ \text{squarefree}}} \frac{L(s, \pi \otimes \chi_D)}{|D|^w} \chi_r(D).
\end{equation}

By the decomposition of the local factors given in Lemma \ref{lem:epsilon-factor-relation-DE} and the constancy of $\chi_{D,\nu} = \chi_{E, \nu}$ for ideals $D \in [E]$ and $\nu \in S$, we have
\begin{equation}
L(s, \pi \otimes \chi_D) = L_S(s, \pi \otimes \chi_E) L^S(s, \pi \otimes \chi_D),
\end{equation}

so that we can rewrite
\begin{equation}
Z_{\star}(s, w; \pi, 1, \delta_E \chi_r) = L_S(s, \pi \otimes \chi_E) \sum_{\substack{D \in [E] \\ \text{squarefree}}} \frac{L^S(s, \pi \otimes \chi_D)}{|D|^w} \chi_r(D).
\end{equation}

We recognize the series on the right as $Z_\star^S(s, w; \pi,1,\delta_E \chi_r) $. The above justifies that we can first concentrate on studying it instead of $Z_\star(s, w; \pi, 1, \delta_E \chi_r) $.

\begin{prop}
\label{prop:meromorphic-continuation-Zstar0}
The function $\xi(s,w) Z_\star^S(s, w; \pi, 1, \delta_E \chi_r)$ admits an analytic continuation  in the region
\begin{equation}
\left\{ \sigma > \tfrac{1}{2}, \tau >  \tfrac{119}{124} \right\} \bigcup \left\{0 \leqslant \sigma \leqslant \tfrac{1}{2}, \tau >  -\tfrac{191}{62} \sigma + \tfrac{5}{2} \right\}.
\end{equation}
\end{prop}

\proof $\bullet$ Let us begin with the case $r = \mathcal{O}$, so that
\begin{equation}
Z_\star^S(s, w; \pi, 1, \delta_E) = \sum_{\substack{D \in [E] \\ \text{squarefree}}} \frac{L^S(s, \pi \otimes \chi_D)}{|D|^w}.
\end{equation}

By the orthogonality relations in $\widehat{H}_C$, it is a uniformly finite linear combination of functions of type $Z^S_\star(s, w; \pi, \rho, 1)$ for $\rho \in \widehat{H}_C$, so that by Proposition \ref{prop:meromorphic-continuation-Zstar}, it also admits a meromorphic continuation to the region \eqref{region}.

$\bullet$ Let us turn to the case where $r$ is a prime ideal. Writing $S_r$ for $S \cup \{r\}$, we get
\begin{align}
Z_\star^S(s, w; \pi, 1, \delta_E \chi_r) &  =  \sum_{\substack{D \in [E] \\ \text{squarefree}}} \frac{L^{S}(s, \pi \otimes \chi_D)}{|D|^w} \chi_r(D) \notag \\
& \label{www} =  \sum_{\substack{D \in [E] \\ \text{squarefree}}} \frac{L^{S_r}(s, \pi \otimes \chi_D)}{|D|^w} \chi_r(D) L_r(s, \pi \otimes \chi_D).
\end{align}

Recalling that $\chi_r$ is a quadratic character so that it only takes the two values $1$ and $-1$ according to the parity of the power of $r$ evaluated, we have 
\begin{equation}
L_r(s, \pi \otimes \chi_D) = \sum_{k \geqslant 0} \frac{a_\pi \left(r^k\right)}{\left|r^k\right|^s} \chi_D\left(r^k\right) = \chi_D(r)L_{1, r}(s) + L_{2, r}(s),
\end{equation}

\noindent where we defined the two partial series according to parity
\begin{equation}
L_{1, r}(s) = \sum_{k \geqslant 0} \frac{a_\pi \left(r^{2k+1}\right)}{\left|r^{2k+1}\right|^s} \qquad \text{and} \qquad L_{2, r}(s) = \sum_{k \geqslant 0} \frac{a_\pi \left(r^{2k}\right)}{\left|r^{2k}\right|^s}.
\end{equation}

By the Blomer-Brumley bounds \cite{blomer_ramanujan_2011} for the Fourier coefficients of automorphic forms on $\GL(3)$, we have that $a_\pi(N) \ll |N|^{5/14+\varepsilon}$ so that both series $L_{1, r}(s) $ and $L_{2, r}(s)$ absolutely uniformly converge on the half-plane $\sigma > 5/14 + \varepsilon$.

By the quadratic reciprocity law \eqref{quadratic-reciprocity-law}, since the sum is over a fixed class in $H_C$ with representative $E$, we have, for all $D \in [E]$,
\begin{equation}
\chi_D(r) = \chi_r(D) \eta(E, r).
\end{equation}

Thus, \eqref{www} rewrites
\begin{align*}
Z_\star^S(s, w; \pi, 1, \delta_E \chi_r) & = \eta(E, r) L_{1, r}(s) \sum_{\substack{D \in [E] \\ (D, r)=1 \\ \text{squarefree}}} \frac{L^{S_r}(s, \pi \otimes \chi_D)}{|D|^w} +  L_{2, r}(s) \sum_{\substack{D \in [E] \\ (D, r)=1 \\ \text{squarefree}}} \frac{L^{S_r}(s, \pi \otimes \chi_D)}{|D|^w} \chi_r(D)
\\
& = \eta(E, r) L_{1, r}(s) Z^{S_r}_\star(s, w; \pi, 1, \delta_E) + L_{2, r}(s) Z^{S_r}_\star(s, w; \pi, 1, \delta_E \chi_r).
\end{align*}

By the analytic continuation of the functions $Z^{S}_{\star}(s, w; \pi, \alpha, \beta)$ around $(\frac{1}{2}, 1)$ obtained in Proposition \ref{prop:meromorphic-continuation-Zstar} we conclude that $Z_\star^S(s, w; \pi, 1, \delta_E \chi_r)$ admits a meromorphic continuation around $(\frac{1}{2}, 1)$, and more precisely that the same double Dirichlet series completed by the factor $\xi(s,w)$ admits analytic continuation on the desired region. \qed

\subsection{Computation of residues at $w=1$}
\label{subsec:residues-computation}
\label{sec:residues}

We aim at understanding the residue
\begin{equation}
\underset{w=1}{\mathrm{Res}} \ Z_\star^S(s, w; \pi, 1,  \delta_E \chi_r) = \lim_{w \to 1} (w-1) Z_\star^S(s, w; \pi, 1, \delta_E \chi_r).
\end{equation}

%

\begin{prop}
We have
\begin{equation}
Z_\star^S(s, w; \pi, 1, \delta_E \chi_r) = h_C^{-1} \sum_{\rho \in \widehat{H}_C} \rho^{-1}(E) \sum_{N \in I(S)} \frac{a_\pi(N)}{|N|^s}  \eta(E, N_0) L(w, \rho\chi_r\chi_{N_0}).
\end{equation}
\end{prop}

\proof The series $Z_\star^S(s, w; \pi, 1, \delta_E \chi_r)$ absolutely converges for $\sigma, \tau > 1$. In this region, we can therefore interchange the summation and recognize $\GL(1)$ $L$-functions for which residues are explicitly computable. More precisely, when $r$ is a prime ideal in $I^+(S)$, 
\begin{align*}
Z_\star^S(s, w; \pi, 1, \delta_E \chi_r) & = \sum_{\substack{D \in [E] \\ \text{squarefree}}} \frac{L^S(s, \pi \otimes \chi_D)}{|D|^w}\chi_r(D) \\
& = \sum_{\substack{D \in [E] \\ \text{squarefree}}} \frac{\chi_r(D)}{|D|^w} \sum_{N \in I^+(S)} \frac{a_\pi(N) \chi_D(N)}{|N|^s} \\
& = \sum_{N} \frac{a_\pi(N)}{|N|^s} \sum_{\substack{D \in [E] \\ \text{squarefree}}} \frac{\chi_r(D)\chi_D(N)}{|D|^w}.
\end{align*}

Formally, we want to apply the quadratic reciprocity law to recognize the rightmost series as a $\GL(1)$ $L$-function. More precisely, decompose $N = N_0 N_1^2 N_2^2$ with $N_0$ the squarefree part of $N$ and $N_1$ such that $P | N_1 \Longrightarrow P | N_0$. In other words, $N_1$ is the square part of $N$ with the same prime factors than $N_0$. If $(D, N_2) \neq 1$, since $D$ is squarefree, we have $\chi_D(N_2)=0$ and therefore $\chi_D(N) = \chi_D(N_2)^2\chi_D(N_0 N_1^2)=0$. We thus assume from now on that $(D, N_2)=1$. If $(D, N_0)=1$ we have, by the quadratic reciprocity law, 
\begin{equation}
\chi_D(N) = \chi_D(N_0) = \eta(E, N_0) \chi_{N_0}(D).
\end{equation}

If $(D, N_0) \neq 1$, we have
\begin{equation}
\chi_D(N) = \chi_D(N_0) = 0 = \eta(E, N_0) \chi_{N_0}(D),
\end{equation}

\noindent so the relation remains true. We can therefore write, using the orthogonality relations, 
\begin{align*}
Z_\star^S(s, w; \pi, 1, \delta_E \chi_r) & = \sum_N \frac{a_\pi(N)}{|N|^s} \sum_{\substack{D \in [E] \\ (D, N) = 1 \\ \text{squarefree}}} \eta(E, N_0) \frac{\chi_r(D) \chi_{N_0}(D)}{|D|^w} \\
& = \sum_N \frac{a_\pi(N)}{|N|^s} h_C^{-1} \sum_{\rho \in \widehat{H}_C} \rho^{-1}(E) \eta(E, N_0) \sum_{\substack{D \\ (D, N) = 1 \\ \text{squarefree}}} \frac{\rho(D)\chi_r(D) \chi_{N_0}(D)}{|D|^w}.
\end{align*}

We recognize the the innermost sum is the $\GL(1)$ $L$-function $L(w, \rho \chi_{r} \chi_{N_0})$. \qed

We can now relate the residue we are interested in in terms of residues of $\GL(1)$ $L$-functions and of the Dedekind zeta function associated to the number field $F$.

\begin{lem}
\label{lem:technica$L$-switching}
We have
\begin{equation}
\begin{split}
& (1+|r|^{-1})^{-1} \sum_{(N_2, r)=1} \frac{a_\pi(N_2^2)}{|N_2^2|^s} \prod_{P | N_2} \left( 1 + \frac{1}{|P|} \right)^{-1} = \frac{1}{|r|^{-1}+ L_{2, r}(s)}  \sum_N  \frac{a_\pi(N^2)}{|N^2|^s}  \prod_{P | N} \left( 1 + \frac{1}{|P|} \right)^{-1}.
\end{split}
\end{equation}
\end{lem}

\proof This is exactly \cite[Equation (4.7)]{chinta_determination_2005}. \qed

\begin{prop}
\label{prop:residues}
Introduce $A(F)$ the residue at $1$ of the Dedekind zeta function $\zeta_F$ associated to $F$. For $\rho \in \widehat{H}_C$, we have,
\begin{equation}
\lim_{w \to 1} (w-1) \zeta_F(2w) \sum_{\substack{D \\ (D, N) = 1 \\ \text{squarefree}}} \frac{\rho \chi_{rN_0}(D)}{|D|^w} = \left\{
\begin{array}{cl}
\displaystyle \prod_{\substack{p \in S \\ \text{or } p | N}} \left( 1  + \frac{1}{|P|^{w}} \right)^{-1} A(F)  & \text{if} \quad \rho \chi_{rN_0} = 1 ; \\
0 & \text{otherwise.} 
\end{array}
\right.
\end{equation}
\end{prop}

\proof Recall that the Dedekind zeta function attached to $F$ is defined by
\begin{equation}
\zeta_F(w) = \prod_p \left( 1 - \frac{1}{|P|^w} \right)^{-1}.
\end{equation}

Moreover, we have the Euler product expansion in the domain of convergence $\tau > 1$, 
\begin{equation}
\sum_{\substack{D \\ (D, N) = 1 \\ \text{squarefree}}} \frac{\rho \chi_{rN_0}(D)}{|D|^w} = \prod_{\substack{p \notin S \\ p \nmid N}} \left( 1 + \frac{\rho \chi_{rN_0}(P)}{|P|^w} \right).
\end{equation}

We therefore have the relation
\begin{equation}
\zeta_F(2w) \sum_{\substack{D \\ (D, N_2) = 1 \\ \text{squarefree}}} \frac{\rho \chi_{rN_0}(D)}{|D|^w}  = \prod_{\substack{p \in S \\ \text{or }P | N}} \left( 1 - \frac{1}{|P|^{2w}} \right)^{-1} \prod_{\substack{p \notin S \\  p \nmid N}} \left( 1  -\frac{\rho \chi_{rN_0}(P)}{|P|^{w}} \right)^{-1}.
\end{equation}

For $\rho \chi_{rN_0} \neq 1$, we recognize the $\GL(1)$ partial $L$-series
\begin{equation}
\prod_{\substack{p \notin S \\ \text{or } p \nmid N}} \left( 1  -\frac{\rho \chi_{rN_0}(P)}{|P|^{w}} \right)^{-1} = L^{S_{N}}(w, \rho \chi_{rN_0}),
\end{equation}

\noindent which is holomorphic on $\tau > 0$. As such, there is no pole on this half-plane containing $w=1$, and we deduce
\begin{equation}
\lim_{w \to 1} (w-1) \zeta_F(2w) \sum_{\substack{D \\ (D, N) = 1 \\ \text{squarefree}}} \frac{\rho \chi_{rN_0}(D)}{|D|^w} = 0.
\end{equation}

For $\rho \chi_{rN_0}= 1$, we have $N_0 = r$. We can thus write $N = r^{2k+1} N_2^2$ for a certain $k \geqslant 0$ and $N_1 \in I^+(S)$ prime to $r$. Explicitly we can therefore compute
\begin{equation}
\zeta_F(2w) \sum_{\substack{D \\ (D, N) = 1 \\ \text{squarefree}}} \frac{\rho \chi_{rN_0}(D)}{|D|^w} = \prod_{\substack{p \in S \\ \text{or } p | N}} \left( 1  + \frac{1}{|P|^{w}} \right)^{-1} \zeta_F(w), 
\end{equation}

\noindent so that, 
\begin{equation}
\lim_{w \to 1} (w-1) \zeta_F(2w) \sum_{\substack{D \\ (D, N) = 1 \\ \text{squarefree}}} \frac{\rho \chi_{rN_0}(D)}{|D|^w} = \prod_{\substack{p \in S \\ \text{or } p | N}} \left( 1  + \frac{1}{|P|^{w}} \right)^{-1} A(F).
\end{equation}

\noindent This proves the claim. \qed

\begin{prop}
\label{prop:residue-explicit-formula}
We have
\begin{align*}
& \underset{w = 1}{\mathrm{Res}} \ Z_\star^S(s, w; \pi, 1, \delta_E \chi_r)  \\
&  = \frac{A(F)}{h_C} \frac{\eta(E, r)}{\zeta_F(2)} \prod_{p \in S} \left(1+\frac{1}{|P|} \right)^{-1} \sum_N \frac{a_\pi(N^2)}{|N^2|^s} \prod_{P | N} \left(1+\frac{1}{|P|} \right)^{-1} (1+|r|^{-1}) \frac{L_{1, r}(s)}{\frac{1}{|r|} + L_{2, r}(s)}.
\end{align*}
\end{prop}

\proof Coming back to $Z_\star^S(s, w; \pi, \alpha, \beta)$ we have
\begin{align*}
& \lim_{w \to 1} \zeta_F(2w) (w-1) Z_\star^S(s, w; \pi, 1, \delta_E \chi_r)  = \\
& \qquad  h_C^{-1} \sum_{\rho \in \widehat{H}_C} \rho^{-1}(E)\sum_N \frac{a_\pi(N)}{|N|^s} \eta(E, N_0)  \lim_{w\to 1}  (w-1) \zeta_F(2w) \sum_{\substack{D \\ (D, N) = 1 \\ \text{squarefree}}} \frac{\rho \chi_r \chi_{N_0}(D)}{|D|^w} .
\end{align*}

It is possible to switch the limit and the summation for $\sigma$ large enough so that we get uniform convergence of the sum over $N$. Finally we get, by Proposition \ref{prop:residues},


\begin{align*}
& \lim_{w \to 1} \zeta_F(2w) (w-1) Z_\star^S(s, w; \pi, 1, \delta_E \chi_r)  \\
& \qquad = h_C^{-1} \sum_{\rho \in \widehat{H}_C} \rho^{-1}(E) \sum_N \frac{a_\pi(N)}{|N|^s} \eta(E, N_0) \lim_{w\to 1} (w-1) \zeta_F(2w)  \sum_{\substack{(D, N) = 1 \\ \text{squarefree}}} \frac{(\rho \chi_r \chi_{N_0})(D)}{|D|^w} \\
& \qquad = h_C^{-1} \eta(E, r) \sum_{\substack{N_2 \\ k \geqslant 0 \\ (N_2, r)=1}} \frac{a_\pi(r^{2k+1}N_2^2)}{|r^{2k+1}N_2^2|^s} \prod_{p \in S} \left( 1 + \frac{1}{|P|} \right)^{-1} \prod_{P | N_2} \left( 1 + \frac{1}{|P|} \right)^{-1} A(F) \\
& \qquad = h_C^{-1}A(F) \eta(E, r) \prod_{p \in S} \left(1+\frac{1}{|P|} \right)^{-1} L_{1, r}(s) \sum_{N_2} \frac{a_\pi(N_2^2)}{|N_2^2|^s} \prod_{P | N_2} \left(1+\frac{1}{|P|} \right)^{-1}.
\end{align*}

Finally we get, dividing by the known residues above, 
\begin{align*}
& \underset{w = 1}{\mathrm{Res}} \ Z_\star^S(s, w; \pi, 1, \delta_E \chi_r)  \\
& = \frac{A(F)}{h_C} \frac{\eta(E, r)}{\zeta_F(2)} \prod_{p \in S} \left(1+\frac{1}{|P|} \right)^{-1} \sum_N \frac{a_\pi(N^2)}{|N^2|^s} \prod_{P | N} \left(1+\frac{1}{|P|} \right)^{-1} (1+|r|^{-1}) \frac{L_{1, r}(s)}{\frac{1}{|r|} + L_{2, r}(s)}.
\end{align*}

This is an explicit expression of the residue in terms of places in $S$. \qed

In order to handle the computations to come, we need a technical results concerning the quantities appearing above.

\begin{lem}
\label{lem:L9}
We have 
\begin{equation}
\sum_N \frac{a_\pi(N^2)}{|N^2|^s} \prod_{P | N} \left(1 + \frac{1}{|P|} \right)^{-1} = L^S(2s, \pi, \mathrm{sym}^2) T(s),
\end{equation}
where $T(s)$ is an absolutely convergent Euler product with no pole around $s=\tfrac{1}{2}$.
\end{lem}

\proof The formal equality is deduced from the computations in the proof of \cite[Proposition 4.2]{chinta_determination_2005}. In the case of Gelbart-Jacquet lifts, we use the bound $a_\pi(N) \ll |N|^{5/14}$ on $\mathrm{GL}(3)$ coefficients and get convergence of the expression and analytic continuation for $\sigma > 3/7$. In the case of non-Gelbart-Jacquet lifts, the convergence and absence of pole is exactly Hypothesis \ref{hyp}. \qed

\textit{Remark.} This is the precise place where we need Hypothesis \ref{hyp}.

\begin{lem}
\label{lem:explicit-functions-L1L2incoeff}
The functions $L_{1, r}(s)$ and $L_{2, r}(s)$ are rational fractions of the Satake parameters $\gamma_{1}(P)$, $\gamma_{2}(P)$, $\gamma_3(P)$ and the Fourier coefficients and $a_\pi(r)$ of $\pi$. More precisely, 
\begin{align*}
L_{1, r}(s) & = \frac{a_\pi(r) + r^{-2s}}{r^s} \prod_{j=1}^3 \left( 1 - \frac{\gamma_j(r)^2}{|r|^{2s}} \right)^{-1} \\
L_{2, r}(s) & = (1 + a_\pi(r)r^{-2s}) \prod_{j=1}^3 \left( 1 - \frac{\gamma_j(r)^2}{|r|^{2s}} \right)^{-1} .
\end{align*}
Moreover, the expression $(1+|r|^{-1}) \frac{L_{1, r}(s)}{|r|^{-1} + L_{2, r}(s)}$ is monotone for $|r|$ sufficiently large and $|\Im s| < |r|^{1-\varepsilon}$.
\end{lem}

\proof This is \cite[Lemma 2.3 and Lemma 5.1]{chinta_determination_2005}. \qed

Coming back to the pure double Dirichlet series $Z_\star(s, w; \pi, 1, \delta_E \chi_r) $, recall that taking aside the places of $S$ we can write
\begin{equation}
Z_\star(s, w; \pi, 1, \delta_E \chi_r) = L_S(s, \pi \otimes \chi_E) Z_\star^S(s, w,\pi, 1, \delta_E \chi_r),
\end{equation}

\noindent and taking the residues at $w=1$, we get
\begin{equation}
\underset{w = 1}{\mathrm{Res}} \ Z_\star(s, w; \pi, 1, \delta_E \chi_r) =  L_S(s, \pi \otimes \chi_E) \cdot  \underset{w = 1}{\mathrm{Res}} \ Z_\star^S(s, w,\pi, 1, \delta_E \chi_r).
\end{equation}

We can summarize the explicit dependencies on $r$ of these residues.

\begin{prop}
\label{prop:residues-Zstar}
We have, for every prime ideal $r$,
\begin{align*}
& \underset{w = 1}{\mathrm{Res}} \ Z_\star(s, w; \pi, 1, \delta_E) = R_1(s, \pi, E) , \\
& \underset{w = 1}{\mathrm{Res}} \ Z_\star(s, w; \pi, 1, \delta_E \chi_r) = R_1(s, \pi, E) R_r(s, \pi, E) .
\end{align*}

\noindent where, for $T(s)$ being the function defined in Lemma \ref{lem:L9}, 
\begin{align}
\label{r1} R_1(s, \pi, E) & = \frac{A(F)}{h_C\zeta_F(2)} \prod_{p \in S} \left(1+\frac{1}{|P|} \right)^{-1} L_S(s, \pi \otimes \chi_E) L^S(2s, \pi, \mathrm{sym}^2) T(s),  \\
\label{rr} R_r(s, \pi, E) & =\eta(E, r)(1+|r|^{-1})\frac{L_{1, r}(s)}{\frac{1}{|r|} + L_{2, r}(s)}.
\end{align}
\end{prop}

\proof Proposition \ref{prop:residue-explicit-formula} states the explicit formula
\begin{align*}
& \underset{w = 1}{\mathrm{Res}} \ Z_\star^S(s, w; \pi, 1, \delta_E \chi_r) =\frac{A(F)}{h_C} \frac{\eta(E, r)}{\zeta_F(2)} \prod_{p \in S} \left(1+\frac{1}{|P|} \right)^{-1} \sum_N \frac{a_\pi(N^2)}{|N^2|^s} \prod_{P | N} \left(1+\frac{1}{|P|} \right)^{-1} (1+|r|^{-1}) \frac{L_{1, r}(s)}{\frac{1}{|r|} + L_{2, r}(s)},
\end{align*}

\noindent so the proposition is a straightforward application of Lemmata \ref{lem:technica$L$-switching} and \ref{lem:L9}. \qed

\textit{Remark.} The partial $L$-factor $L_S(s, \pi \otimes \chi_E)$ is omitted in \cite{chinta_determination_2005}. Even though it is present in their expression (4.9), it should be added back from (4.12) until the end of the paper. In particular, it modifies the definition of their $C_M(\pi)$ in (4.21). This is not a critical issue since the partial $L$-factor $L_S(s, \pi \otimes \chi_E)$ never vanishes by the Euler product expression and the known bounds towards the Ramanujan conjecture. We are grateful to Adrian Diaconu for having clarified this fact to us.

Coming back to the local behavior at $(\frac{1}{2}, 1)$, recall that the analytic continuation stated by Proposition \ref{prop:meromorphic-continuation-Zstar0} implies the local development
\begin{equation}
\label{loca$L$-development}
Z_\star^S(s, w; \pi, 1, \delta_E \chi_r) = \frac{A(s)}{w-1} + \frac{B(s)}{w+3s-5/2} +H(s, w),
\end{equation}

\noindent where $H(s, w)$ is a holomorphic function around the point $(\frac{1}{2}, 1)$.

\subsection{Local analysis}
\label{subsec:fourier-coefficients}


\subsubsection{Case of Gelbart-Jacquet lifts}

Suppose now that $\pi$ is an automorphic cuspidal self-contragredient representation of $\GL(3)$ with trivial central character, and that $\pi$ is a Gelbart-Jacquet lift. 

\begin{prop}
\label{prop:GJ}
For $\pi$ a Gelbart-Jacquet lift, there is a ray class $E$ such that
\begin{equation}
\lim_{s\to \frac{1}{2}} \lim_{w  \to 1} (w-1) \left( s - \tfrac{1}{2} \right) Z_\star(s, w; \pi, 1, \delta_E \chi_r) = \left( \underset{s = 1/2}{\mathrm{Res}} R_1\left(s, \pi, E \right) \right)R_r\left(\tfrac{1}{2}, \pi, E\right).
\end{equation}
\end{prop}

\proof By the results of \cite{grs}, $\pi$ being a Gelbart-Jacquet lift is equivalent to $L^S(s, \pi, \mathrm{sym}^2)$ having a pole at $s=1$. Moreover, the partial $L$-factor $L_S(s, \pi \otimes \chi_E)$ never vanishes by the Euler product expression and the known bounds towards the Ramanujan conjecture. Then \eqref{r1} implies that $R_1(s, \pi, E)$ has a pole at $s=1/2$. We can therefore write 
\begin{align*}
A(s) & = \frac{A_0}{s-\frac{1}{2}} + A_1(s),  \\
B(s) & = \frac{B_0}{s-\frac{1}{2}} + B_1(s),
\end{align*}

where $A_1$ and $B_1$ are holomorphic functions around $s=\frac{1}{2}$. This allows for a refinement of the local expansion \eqref{loca$L$-development} into
\begin{equation}
\label{le}
\begin{split}
Z_\star^S(s, w; \pi, 1, \delta_E \chi_r) & = \frac{A_0}{(w-1)(s-\frac{1}{2})} + \frac{A_1(s)}{w-1}  + \frac{B_0}{(w+3s-\frac{5}{2})(s-\frac{1}{2})} + \frac{B_1(s)}{w+3s-\frac{5}{2}} + H(s, w).
\end{split}
\end{equation}

Now, we use the functional equation 
\begin{equation}
Z_\star^S(s, w; \pi, 1, \delta_E \chi_r) = M(s, E) Z_\star^S(\phi(s, w); \pi, 1, \delta_E \chi_r),
\end{equation}
where
\begin{equation}
  M(s,E) = \left[ \frac{\varepsilon(s, \pi \otimes \chi_E)}{|E_0|^{3(1/2-s)}} \right] \prod_{v \in S} \frac{L_v(1-s, \pi \otimes \chi_E)}{L_v(s, \pi\otimes \chi_E)}.
\end{equation}

By \cite[Theorem 0.3 (ii)]{bump_sums_2004}, there are infinitely many character twists such that $L^S(s, \pi \otimes \chi_E)$ does not vanish at $s=1/2$. For such a ray class $E$, we see that $M(1/2, E) = 1$. Applying the functional equation above to \eqref{le} yields 
\begin{align*}
Z_\star^S(s, w; \pi, 1, \delta_E \chi_r) & = \frac{A_0}{(w-1)(s-\frac{1}{2})} + \frac{A_1(s)}{w-1} + \frac{B_0}{(w+3s-\frac{5}{2})(s-\frac{1}{2})} + \frac{B_1(s)}{w+3s-\frac{5}{2}} + H(s, w) \\
& = M(s, E) \left[ - \frac{A_0}{(w+3s-\frac{5}{2})(s-\frac{1}{2})} + \frac{A_1(1-s)}{w+3s - \frac{5}{2}} -  \frac{B_0}{(w-1)(s-\frac{1}{2})} + \frac{B_1(1-s)}{w-1} + H(\phi(s, w)) \right].
\end{align*}

Comparing polar parts give 
\begin{equation*}
A_0 =  -B_0.
\end{equation*}

With this, \eqref{le} simplifies to
\begin{equation}
Z_\star^S(s, w; \pi, 1, \delta_E \chi_r) = \frac{3A_0}{(w-1)(w+3s-\frac{5}{2})} + \frac{A_1(s)}{w-1} + + \frac{B_1(s)}{w+3s-\frac{5}{2}} + H(s, w).
\end{equation} 

Taking the residue at $w=1$ and adding back the local factor $L_S(s, \pi \otimes \chi_E)$, we get
\begin{equation}
 \underset{w = 1}{\mathrm{Res}} \ Z_\star \left(s, w; \pi, 1, \delta_E \chi_r\right)  = \frac{A_0 L_S(s, \pi \otimes \chi_E)}{s-\frac{1}{2}} + A_1(s) L_S(s, \pi \otimes \chi_E).
\end{equation}

Examining the residue at $s=1/2$ and taking into account Proposition \ref{prop:residues-Zstar} lead to
\begin{equation}
A_0 \  L_S(\tfrac{1}{2}, \pi \otimes \chi_E) = \left( \underset{s = 1/2}{\mathrm{Res}} R_1\left(s, \pi, E\right) \right) R_r\left(\tfrac{1}{2}, \pi, E\right)
\end{equation}

\noindent as claimed. \qed

Note in particular that $\underset{s = 1/2}{\mathrm{Res}} R_1(s, \pi, E)$ is nonzero. Indeed, there is a simple pole at $s=1/2$ for $L^S(2s, \pi, \mathrm{sym}^2)$ implying the corresponding residue is nonzero. Moreover, the partial $L$-factor $L_S(s, \pi \otimes \chi_E)$ does not vanish at $s=1/2$ by the Euler product expression and the known bounds towards the Ramanujan conjecture. 

\subsubsection{Case of non-Gelbart-Jacquet lifts}

Assume $\pi$ is an automorphic cuspidal representation of $\GL(3)$ with trivial central character, and that $\pi$ is not a Gelbart-Jacquet lift. 

\begin{prop}
\label{prop:nGJ}
If $\pi$ is a non-Gelbart-Jacquet lift and there is a ray class $E$ such that $L(1/2, \pi \otimes \chi_E) \neq 0$, then
\begin{equation}
\lim_{w  \to 1} (w-1)  Z_\star\left(\tfrac{1}{2}, w; \pi, 1,  \delta_E \chi_r\right) = 2R_1\left(\tfrac{1}{2}, \pi, E\right)R_r\left(\tfrac{1}{2}, \pi, E\right).
\end{equation}
\end{prop}

\proof By the results of \cite{grs}, the condition on $\pi$ is equivalent to $L^S(s, \pi, \mathrm{sym}^2)$ having no pole at $s=1$. Since the partial $L$-factor $L_S(s, \pi \otimes \chi_E)$ is a finite Euler product, it also has no pole at $s=1/2$. Then \eqref{r1} implies that $R_1(s, \pi, E)$ has no pole at $s=1/2$.  

Now, we use the functional equation 
\begin{equation}
Z_\star^S(s, w; \pi, 1, \delta_E \chi_r) = M(s,E) Z_\star^S(\phi(s, w); \pi, 1, \delta_E \chi_r) .
\end{equation}

By the non-vanishing assumption, we deduce that $M(1/2, E)=1$. Applying the functional equation above to \eqref{loca$L$-development} therefore yields
\begin{align*}
Z_\star^S(s, w; \pi, 1, \delta_E \chi_r) & = \frac{A(s)}{w-1} + \frac{B(s)}{w+3s-5/2} + H(s, w)\\
& = M(s, E) \left[ \frac{A(1-s)}{w + 3s-5/2} + \frac{B(1-s)}{w-1} + H(\phi(s, w)) \right],
\end{align*}

so that in particular we deduce
\begin{equation}
B(s) = M(s, E)  A(1-s) \text{ and } A(1/2) = B(1/2).
\end{equation}

 Letting $s \to 1/2$ in \eqref{loca$L$-development} yields
\begin{equation}
Z_\star\left( \frac{1}{2}, w; \pi, 1,  \delta_E \chi_r \right) = \frac{2 R_{1}\left( \frac{1}{2}, \pi, E \right)R_r\left( \frac{1}{2}, \pi, E\right)}{w-1} + H\left( \frac{1}{2}, w \right) L_S \left( \frac{1}{2}, \pi \otimes \chi_E \right),
\end{equation}

so that taking the residue at $w=1$ yields the result. \qed

We will need to know that $R_1(s, \pi, E)$ does not vanish at $s=\frac{1}{2}$ for the final argument. We summarize sufficient conditions for this in the following lemma.
\begin{lem}
\label{lem:R1-nonzero}
If $\pi$ is not a Gelbart-Jacquet lift and if there is a class $E$ such that $L(1/2, \pi \otimes \chi_E) \neq 0$, then $R_1(1/2, \pi, E) \neq 0$. 
\end{lem}

\proof The symmetric square L-function $L(s, \pi, \mathrm{sym}^2)$ does not vanish at $s=1$, see e.g. \cite[Theorem A.1]{lapid_harish-chandra_2012}. The result is now a straightforward consequence of \eqref{r1}.  \qed

\subsection{Proof of Theorem \ref{thm:result-nv}}
\label{nonvanishing}

Let $\pi$ be a cuspidal automorphic representation of $\GL(3)$ that is not a Gelbart-Jacquet lift. In particular, that implies that $L(s, \pi, \mathrm{sym}^2)$ has no pole at $s=1$. For all $w \neq 1$ near $1$, we can take $s=1/2$ and recall the local expansion proven above: there is an analytic function $H(w)$ around $1$ such that
\begin{equation}
Z_\star^S \left( \frac{1}{2}, w ; \pi, 1, \delta_E \right) = \frac{2 R_1\left( \frac{1}{2}, \pi, E \right)}{L_S(1/2, \pi \otimes \chi_E)(w-1)} + H\left( \frac{1}{2}, w \right).
\end{equation}

Assume that there is a class $E$ such that $L^S(1/2, \pi \otimes \chi_E) \neq 0$. By Lemma \ref{lem:R1-nonzero}, $R_1(1/2, \pi, E)$ does not vanish. In particular, $Z^S( 1/2, w ; \pi, \delta_E, 1 )$ has a pole at $w=1$.

We now consider the following smoothed truncated version of the double Dirichlet series. For any $x>0$, let
\begin{equation}
I(x) = \sum_{\substack{D \in I^+(S) \\ D \text{ squarefree}}} L^S\left( \frac{1}{2}, \pi \otimes \chi_D \right) e^{-|D|/x}.
\end{equation}

Recall the expression, following \cite[Theorem 3.8]{bump_sums_2004},
\begin{equation}
e^{-1/x} = \frac{1}{2\pi i}\int_{2-i\infty}^{2+i\infty} \Gamma(x) x^w \mathrm{d}w, 
\end{equation}

which allows for an integral formulation of $I(x)$ given by
\begin{equation}
I(x) = \frac{1}{2\pi i} \int_{2-i\infty}^{2+i\infty} Z_\star^S(1/2, w; \pi, 1, \delta_E) \Gamma(w) x^w \mathrm{d}w.
\end{equation}

We can therefore move the vertical integration line to $\Re(w) = 119/124 +  \varepsilon$ where the integral converges absolutely and contributes an error term of order $x^{119/124+\varepsilon}$ by the bounds on $\Gamma$ and the vertical estimates provided in Section \ref{subsec:vertica$L$-bounds}. In this process, we pick up the residue corresponding to the pole of $Z_\star^S(1/2, w ; \pi, 1, \delta_E)$ at $w=1$. The residue at $w=1$ contributes
\begin{equation}
\underset{w=1}{\mathrm{Res}} \ Z_\star^S(1/2, w ; \pi, 1, \delta_E) \Gamma(w) x^w = \frac{2 R_1\left( \frac{1}{2}, \pi, E \right)}{L_S(1/2, \pi \otimes \chi_E)} x.
\end{equation}
We know from Lemma \ref{lem:R1-nonzero} that $R_1(1/2,\pi,E) \neq 0$. We note that the contribution of a squarefree integral ideal $D \in I^+(S)$ with $|D| \leqslant x$ is of size
\begin{align*}
L^S\left( \tfrac{1}{2}, \pi \otimes \chi_D \right) & \ll |D|^{3/4+\varepsilon} \ll x^{3/4 + \varepsilon}.
\end{align*}

Therefore, a finite number of such $L$-factors $L^S(1/2, \pi \otimes \chi_D)$ with $|D| \leqslant x$ cannot alone contribute to a size of $x$. We deduce that there should be infinitely many non-vanishing of the $L$-factors $L^S\left( 1/2, \pi \otimes \chi_D \right)$.

\subsection{Proof of Theorem \ref{thm:result}}
\label{fourier-coeff-proof-theorem}

Let $\pi$ and $\pi'$ be two automorphic representations as in the statement of the theorem, and $E$ an ideal ray class as in the assumptions of Theorem \ref{thm:result}. Let $S$ be large enough in order to the local components of $\pi$ and $\pi'$ to be both unramified and principal series at the places outside $S$. Suppose moreover that there is a nonzero constant $\kappa$ such that
\begin{equation}
L\left(\frac{1}{2}, \pi \otimes \chi_D\right) = \kappa \cdot L\left(\frac{1}{2}, \pi' \otimes \chi_D\right),
\end{equation}

\noindent for every $D \in I^+(S)$. By summing over squarefree $D$ in the same class as $E \in H_C$, we deduce the relation between the associated double Dirichlet series in the region of convergence, for every $r$ prime or equal to 
$\mathcal{O}$,
\begin{equation}
Z_\star\left(\frac{1}{2}, w; \pi, 1, \delta_E \chi_r \right) = \kappa \cdot Z_\star\left(\frac{1}{2}, w; \pi', 1, \delta_E \chi_r\right).
\end{equation}

By Proposition \ref{prop:GJ} or \ref{prop:nGJ} depending on whether or not $\pi$ is a Gelbart-Jacquet lift we get, for $r= \mathcal{O}$, the relevant $R_1$-related qualities are equal and nonzero. It follows from the above that
%

\begin{equation}
R_r\left(\frac{1}{2}, \pi, E\right) = R_r\left(\frac{1}{2}, \pi', E\right).
\end{equation}

By Lemma \ref{lem:explicit-functions-L1L2incoeff}, these quantities are monotone for large enough $|r|$. In particular they are injective and the above equality yields that $\pi$ and $\pi'$ have same Fourier coefficients for large enough $|r|$. By the strong multiplicity one theorem, we conclude that $\pi \cong \pi'$. \qed

\bibliographystyle{apalike}
\bibliography{biblio_twists}

\begin{thebibliography}{}

\bibitem[Blomer and Brumley, 2011]{blomer_ramanujan_2011}
Blomer, V. and Brumley, F. (2011).
\newblock On the {Ramanujan} conjecture over number fields.
\newblock {\em Annals of Mathematics}, 174(1):581--605.
\newblock 1.

\bibitem[Bump, 1997]{bump_automorphic_1997}
Bump, D. (1997).
\newblock {\em Automorphic {Forms} and {Representations}}.
\newblock Number~55 in Cambridge {Studies} in {Adv}. {Math}. Cambridge
  University Press, Cambridge.

\bibitem[Bump et~al., 2004]{bump_sums_2004}
Bump, D., Friedberg, S., and Hoffstein, J. (2004).
\newblock Sums of twisted {GL}(3) automorphic {L}-functions.
\newblock In {\em H. {Hida}, {D}. {Ramakrishnan}, {F}. {Shahidi} ({Eds}.),
  {Contributions} to {Automorphic} {Forms}, {Geometry} and {Arithmetic}}, pages
  131--162. Johns Hopkins University Press.
\newblock 1.

\bibitem[Cassels and Fröhlich, 1967]{cassels_algebraic_1967}
Cassels and Fröhlich (1967).
\newblock {\em Algebraic {Number} {Theory}}.
\newblock Academic Press Inc.

\bibitem[Chinta and Diaconu, 2005]{chinta_determination_2005}
Chinta, G. and Diaconu, A. (2005).
\newblock Determination of a {GL}(3) cuspform by twists of central {L}-values.
\newblock {\em IMRN}, 4(48):1995--2040.

\bibitem[Chinta et~al., 2005]{chinta_asymptotics_2005}
Chinta, G., Friedberg, S., and Hoffstein, J. (2005).
\newblock Asymptotics for sums of twisted {L}-functions and applications.
\newblock In {\em Automorphic {Representations}, {L}-{Functions} and
  {Applications}: {Progress} and {Prospects}}. DE GRUYTER, Berlin, New York.

\bibitem[Cogdell, 2004]{cogdell_analytic_2004}
Cogdell, J.~W. (2004).
\newblock Analytic {Theory} of {L}-functions for {GL}(n).
\newblock In Bernstein, J. and Gelbart, S., editors, {\em An {Introduction} to
  the {Langlands} {Program}}, pages 251--268. Birkhäuser Boston, Boston, MA.

\bibitem[Diaconu et~al., 2003]{diaconu_multiple_2003}
Diaconu, A., Goldfeld, D., and Hoffstein, J. (2003).
\newblock Multiple {Dirichlet} {Series} and {Moments} of {Zeta} and
  {L}-{Functions}.
\newblock {\em Compo. Math.}, 139(3):297--360.

\bibitem[Farmer et~al., 2017]{farmer_analytic_2017}
Farmer, D., Pitale, A., Ryan, N., and Schmidt, R. (2017).
\newblock Analytic {L}-functions: {Definitions}, {Theorems} and {Connections}.
\newblock Preprint available at arXiv:1711.10375.

\bibitem[Fisher and Friedberg, 2003]{fisher_sums_2003}
Fisher, B. and Friedberg, S. (2003).
\newblock Sums of twisted {GL}(2) {L} -functions over function fields.
\newblock {\em Duke Mathematical Journal}, 117(3):543--570.

\bibitem[Fisher and Friedberg, 2004]{fisher_double_2004}
Fisher, B. and Friedberg, S. (2004).
\newblock Double {Dirichlet} series over function fields.
\newblock {\em Compositio Mathematica}, 140(03):613--630.

\bibitem[Friedberg and Hoffstein, 1995]{friedberg_nonvanishing_1995}
Friedberg, S. and Hoffstein, J. (1995).
\newblock Nonvanishing {Theorems} for {Automorphic} {L}-{Functions} on {GL}(2).
\newblock {\em The Annals of Mathematics}, 142(2):385.

\bibitem[Gelbart and Jacquet, 1978]{ASENS_1978_4_11_4_471_0}
Gelbart, S. and Jacquet, H. (1978).
\newblock A relation between automorphic representations of ${\rm gl}(2)$ and
  ${\rm gl}(3)$.
\newblock {\em Annales scientifiques de l'\'Ecole Normale Sup\'erieure}, Ser.
  4, 11(4):471--542.

\bibitem[Gelfand et~al., 1969]{gelfand_representation_1969}
Gelfand, I.~M., Graev, M.~I., and Pyatetskii-Shapiro, I.~I. (1969).
\newblock {\em Representation theory and automorphic functions}.
\newblock Academic Press, Inc.
\newblock Translated from the Russian by K. A. Hirsch. Reprint of the 1969
  edition. Generalized Functions, 6.

\bibitem[Ginzburg et~al., 2005]{cogdell_nonvanishing_2005}
Ginzburg, D., Jiang, D., and Rallis, S. (2005).
\newblock On the nonvanishing of the central value of the {Rankin}-{Selberg}
  {L}-functions {II}.
\newblock In Cogdell, J.~W., Jiang, D., Kudla, S.~S., Soudry, D., and Stanton,
  R.~J., editors, {\em Automorphic {Representations}, {L}-{Functions} and
  {Applications}: {Progress} and {Prospects}}. DE GRUYTER, Berlin, New York.

\bibitem[Ginzburg et~al., 1999]{grs}
Ginzburg, D., Rallis, S., and Soudry, D. (1999).
\newblock On {Explicit} {Lifts} of {Cusp} {Forms} from {GLm} to {Classical}
  {Groups}.
\newblock {\em Annals of Mathematics}, 150(3):807--866.

\bibitem[Goldmakher and Louvel, 2013]{goldmakher_quadratic_2013}
Goldmakher, L. and Louvel, B. (2013).
\newblock A quadratic large sieve inequality over number fields.
\newblock {\em Mathematical Proceedings of the Cambridge Philosophical
  Society}, 154(02):193--212.

\bibitem[Heath-Brown, 1995]{heath-brown_mean_1995}
Heath-Brown, D.~R. (1995).
\newblock A mean value estimate for real character sums.
\newblock {\em Acta Arithmetica}, 72(3):235--275.

\bibitem[Hörmander, 1990]{hormander_introduction_1990}
Hörmander (1990).
\newblock {\em An {Introduction} to {Complex} {Analysis} in {Several}
  {Variables}}.
\newblock North-Holland.

\bibitem[Iwaniec and Sarnak, 2000]{iwaniec_perspectives_2000}
Iwaniec, H. and Sarnak, P. (2000).
\newblock Perspectives on the {Analytic} {Theory} of {L}-functions.
\newblock {\em Geom. Funct. Anal}, pages 705--741.

\bibitem[Kim and Sarnak, 2002]{kim_functoriality_2002}
Kim, H.~H. and Sarnak, P. (2002).
\newblock Functoriality of the exterior square of {GL4} and the symmetric
  fourth of {GL2}.
\newblock {\em Journal of the American Mathematical Society}, 16(1):139--183.

\bibitem[Lapid, 2012]{lapid_harish-chandra_2012}
Lapid, E. (2012).
\newblock On the {Harish}-{Chandra} {Schwartz} space of {G}({F}
  ){\textbackslash}{G}({A}).
\newblock In {\em Proceedings of the {International} {Colloquium} on
  {Automorphic} {Representations} and {L}-{Functions}}, page~36. Mumbai.
\newblock Appendix by F. Brumley, Lower bounds on Rankin-Selberg L-functions.

\bibitem[Li, 2007]{li_determination_2007}
Li, J. (2007).
\newblock Determination of a {GL}2 automorphic cuspidal representation by
  twists of critical {L}-values.
\newblock {\em J. Number Theory}, 123(2):255--289.

\bibitem[Luo, 2005]{luo_nonvanishing_2005}
Luo, W. (2005).
\newblock Nonvanishing of \${L}\$-functions for \${\textbackslash}{GL}(n,
  {\textbackslash}mathbf\{{A}\}\_q)\$.
\newblock {\em Duke Mathematical Journal}, 128(2):199--207.

\bibitem[Luo and Ramakrishnan, 1997]{luo_determination_1997}
Luo, W. and Ramakrishnan, D. (1997).
\newblock Determination of modular forms by twists of critical {L}-values.
\newblock {\em Inventiones mathematicae}, 130(2):371--398.

\bibitem[Montgomery and Vaughan, 2006]{montgomery_multiplicative_2006}
Montgomery, H.~L. and Vaughan, R.~C. (2006).
\newblock {\em Multiplicative {Number} {Theory} {I}: {Classical} {Theory}}.
\newblock Cambridge University Press, Cambridge.

\bibitem[Munshi and Sengupta, 2015]{munshi_determination_2015}
Munshi, R. and Sengupta, J. (2015).
\newblock Determination of {GL}(3) {Hecke}–{Maass} forms from twisted central
  values.
\newblock {\em Journal of Number Theory}, 148:272--287.

\bibitem[Ramakrishnan, 2014]{ramakrishnan_exercise_2014}
Ramakrishnan, D. (2014).
\newblock An exercise concerning the selfdual cusp forms on {GL}(3).
\newblock {\em Indian Journal of Pure and Applied Mathematics}, 45(5):777--785.

\bibitem[Waldspurger, 1985]{waldspurger_sur_1985}
Waldspurger, J.-L. (1985).
\newblock Sur les valeurs de certaines fonctions \${L}\$ automorphes en leur
  centre de symétrie.
\newblock {\em Compositio Mathematica}, 54(2):173--242.

\end{thebibliography}

\end{document}